\documentclass[12pt,a4paper,final]{amsart}
\usepackage{amsmath}
\usepackage{amssymb}
\usepackage{amsthm}
\usepackage{color}
\usepackage{enumerate}
\usepackage[dvips]{graphicx}
\usepackage{comment}

\newcommand{\T}{\mathbb{T}}

\DeclareMathOperator\supp{supp}

\DeclareMathOperator\Lip{Lip}

\DeclareMathOperator\SC{SC}

\DeclareMathOperator\dom{dom}

\def\R{\mathbb{R}}

\def\N{\mathbb{N}}
\def\M{\mathbb{M}}
\def\cF{\mathcal{F}}

\def\bye{\end{document}}
\def\by{\end{proof}\bye}

\def\hello{\begin{document}}
\def\fr{\frac}
\def\disp{\displaystyle}
\def\ga{\alpha}
\def\go{\omega}
\def\gep{\varepsilon}
\def\ep{\gep}
\def\mid{\,:\,}
\def\gb{\beta}
\def\gam{\gamma}
\def\gd{\delta}
\def\gz{\zeta}
\def\gth{\theta}
\def\gk{\kappa}
\def\gl{\lambda}
\def\gL{\Lambda}
\def\gs{\sigma}
\def\gf{\varphi}
\def\tim{\times}
\def\aln{&\,}
\def\ol{\overline}
\def\ul{\underline}
\def\pl{\partial}
\def\hb{\text}
\def\cF{\mathcal{F}}
\def\Int{\mathop{\text{int}}}

\def\gG{\varGamma}
\def\lan{\langle}
\def\ran{\rangle}
\def\cD{\mathcal{D}}
\def\cB{\mathcal{B}}
\def\bcases{\begin{cases}}
\def\ecases{\end{cases}}
\def\balns{\begin{align*}}
\def\ealns{\end{align*}}
\def\balnd{\begin{aligned}}
\def\ealnd{\end{aligned}}
\def\1{\mathbf{1}}
\def\bproof{\begin{proof}}

\def\eproof{\end{proof}}

\theoremstyle{definition}
\newtheorem{definition}{Definition}
\theoremstyle{plain}
\newtheorem{theorem}[definition]{Theorem}
\newtheorem{corollary}[definition]{Corollary}
\newtheorem{lemma}[definition]{Lemma}
\newtheorem{proposition}[definition]{Proposition}
\theoremstyle{remark}
\newtheorem{remark}[definition]{Remark}
\newtheorem{notation}[definition]{Notation}

\def\hr#1{\rule{0pt}{#1 pt}}

\def\red#1{\textcolor{red}{#1}}
\def\blu#1{\textcolor{blue}{#1}}
\def\beq{\begin{equation}}
\def\eeq{\end{equation}}
\def\bthm{\begin{theorem}}
\def\ethm{\end{theorem}}
\def\bproof{\begin{proof}}
\def\eproof{\end{proof}}

\usepackage[abbrev]{amsrefs}
\newcommand\coolrightbrace[2]{%
\left.\vphantom{\begin{matrix} #1 \end{matrix}}\right\}#2}
\def\eqr#1{\eqref{#1}}
\def\bmat{\begin{pmatrix}}
\def\emat{\end{pmatrix}}

\newcommand{\Pmo}{\mathcal{P}^-_{1}}
\newcommand{\Ppo}{\mathcal{P}^+_{1}}
\newcommand{\Pmk}{\mathcal{P}^-_{k}}
\newcommand{\Ppk}{\mathcal{P}^+_{k}}
\def\diag{\operatorname{diag}}
\def\rT{\mathrm{T}}
\newcommand{\Rn}{{\mathbb R}^N}
\def\IN{\text{ in } }

\def\AND{\text{ and }}
\def\FOR{\text{ for }}
\def\FORALL{\text{ for all }}
\def\ON{\text{ on }}
\def\IF{\hb{ if }}
\def\WITH{\text{ with }}

\def\rmb{\mathrm{b}}
\def\I{\mathbb{I}}
\def\du#1{\left\lan#1\right\ran}
\def\bald{\begin{aligned}}
\def\eald{\end{aligned}}
\def\stm{\setminus}
\def\t{\tau}
\def\SC-{\operatorname{SC}^-}
\def\sc-{\operatorname{sc}^-}
\def\lip{\operatorname{lip}}
\def\B{\operatorname{\mathbb{B}}}

\def\Sp{\operatorname{Sp}}
\def\P{\operatorname{\mathbb{P}}}
\def\erf{\eqref}
\def\cG{\mathcal{G}}
\def\cE{\mathcal{E}}
\def\pr{\,^\prime}
\def\gX{\Xi}
\def\gS{\Sigma}
\def\cV{\mathcal{V}}
\def\cW{\mathcal{W}}
\def\Hall{(H_i)_{i\in\I}}\def\Lall{(L_i)_{i\in\I}}
\def\prop{proposition}
\def\thm{theorem}
\def\lem{lemma}

\def\fC{\frak{C}}
\def\fM{\frak{M}}
\def\0{\mathbf{0}}

\renewcommand{\subjclassname}{%
\textup{2010} Mathematics Subject Classification}

\title[Vanishing discount problem]{The vanishing discount problem for \\
monotone systems
of Hamilton-Jacobi equations. \\ Part 2: Nonlinear coupling}

\author[H. Ishii]{Hitoshi Ishii}
\address[\textsc{Hitoshi Ishii}]{Institute for Mathematics and Computer Science\newline
\indent Tsuda University  \newline
 \indent   2-1-1 Tsuda, Kodaira, Tokyo, 187-8577 Japan.
}
\email{hitoshi.ishii@waseda.jp}

\author[L. Jin]{Liang Jin}
\address[\textsc{Liang Jin}]{
\newline
\indent  Department of Mathematics, Nanjing University of Science and Technology  \newline
 \indent  Nanjing 210094,  China.
}
\email{jl@njust.edu.cn}

\keywords{systems of Hamilton-Jacobi equations, Mather measures, vanishing discount}
\subjclass[2010]{
35B40, 
35F50, 
49L25 
}

\def\alert#1{\begin{color}{red}#1 \end{color}}

\begin{document}
\maketitle
\begin{abstract}
We study the vanishing discount problem
for a nonlinear monotone system of Hamilton-Jacobi equations.
This continues the first author's investigation on the vanishing
discount problem for a monotone system of Hamilton-Jacobi equations.
As in Part 1, we introduce by the convex duality Mather measures and their analogues for the system, which we call respectively Mather and Green-Poisson measures, and prove a convergence theorem for the vanishing discount problem.
Moreover, we establish an existence result for the ergodic problem.
\end{abstract}

\tableofcontents

\def\rmc{\mathrm{c}}

\section{Introduction}

We consider the $m$-system of Hamilton-Jacobi equations
\beq\tag{P$_\gl$}\label{Pl}
\gl v_i^\gl(x)+H_i(x,Dv_i^\gl(x),v^\gl(x))=0 \ \ \IN \T^n,\ i\in\I,
\eeq
where $\I:=\{1,\ldots,m\}$ with $m\in\N$, $\gl$ is a nonnegative constant,
called the discount factor in terms of
optimal control.
Here $\T^n$ denotes the $n$-dimensional flat torus and $H=(H_i)_{i\in\I}$ is a family of continuous Hamiltonians. The unknown in \erf{Pl} is an $\R^m$-valued function
$v^\gl=(v^\gl_i)_{i\in \I}$ on $\T^n$ and the above system can be written in the vector form as follows:
\[\tag{P$_\gl$}
\gl v^\gl+H[v^\gl]=0 \ \ \IN \T^n.
\]
We have used here the abbreviated expression $H[v^\gl]$ to denote
$(H_i(x,Dv_i^\gl(x),v^\gl(x)))_{i\in\I}$.
The system is weakly coupled in the sense that every $i$-th equation depends on $Dv^\gl$ only through $Dv_i^\gl$ but not on $Dv_j^\gl$, with $j\not=i$.

We are concerned with the vanishing discount problem for \erf{Pl}, that is,
the asymptotic behavior of the solution $v^\gl$ of
\erf{Pl} as $\gl \to 0+$. Notably, the main concern is the convergence of the whole family $(v^\gl)_{\gl>0}$ as $\gl \to 0+$. 

Recently, there has been a great interest in the vanishing discount problem
concerned with Hamilton-Jacobi equations and, furthermore, fully nonlinear degenerate elliptic PDEs. We refer to \cites{DFIZ, AAIY, MT, IMT1, IMT2, GMT,
CCIZ, DZ2, IS, ItSM, DFIZdis, Go08, CaCDGo} for relevant work.
The asymptotic analysis in these papers relies heavily on Mather measures
or their generalizations and, thus, it is considered part of Aubry-Mather and weak KAM theories. For the development of these theories we refer to
\cites{Fa1, Fa2, Ev1} and the references therein.
We refer to \cite{CGT, CLLN, DSZ, DZ1, Ev2, MSTY, MT1,MT2,MT3} for the recent development in the asymptotic analysis and weak KAM theory
for systems of Hamilton-Jacobi equations.

We are here interested in the case of systems of Hamilton-Jacobi equations.
Davini and Zavidovique in \cite{DZ2} have established
a general convergence result for the vanishing discount problem for \erf{Pl} when
the coupling is linear and the coupling coefficients are constant.
Adapting the convex duality argument in \cite{IMT1}
to the system,  the first author of this paper has treated the case of linear coupling,
with the coupling coefficients depending on the space variable.
In this paper, we extend the scope of the previous work \cite{Icoupling}
and discuss the case of the system with nonlinear
coupling.  Our argument is pretty much parallel to that in \cite{Icoupling}.
We refer for further references to \cite{Icoupling, DZ2}.

Under our hypotheses described later, the 
limit function $v^0$ of the solution $v^\gl$ of \erf{Pl} satisfies the system 
$H[v^0]=0$ and the principal difficulty in the asymptotic analysis lies in the fact that the system $H[u] = 0$, or (P$_0$), usually has multiple solutions of which the structure is not simple in general. The critical role of Mather measures is indeed to identify the limit function $v^0$ from the solutions of $H[u]=0$.

We assume throughout (see \erf{convex} below) that the functions: 
$\R^n\tim\R^m\ni(p,u) \mapsto H_i(x,p,u)$ are convex for $(x,i)\in\T^n\tim\I$.  
We have chosen to make the convexity requirement on $H$ in $u$ simply because 
of our technical limitations for studying the vanishing discount problem.   
The non-convexity issue for the vanishing discount problem 
has already been addressed in \cite{GMT, CCIZ, Zi} 
in the case of scalar equations.  The paper \cite{Zi} illustrates by examples 
that, in general,  the convergence of the whole family $(v^\gl)_{\gl>0}$ does not 
hold without the convexity of $p\mapsto H_i(x,p,u)$, while \cite{GMT, CCIZ} 
indicate some possible generalizations beyond the convexity of 
$u\mapsto H_i(x,p,u)$.    

On the other hand, under the convexity assumption, the system \erf{Pl} 
may be regarded as the dynamic programming equation of optimal control 
of random evolutions, where the state $(x(t),i(t))$ at time $t$ is in $\T^n\tim\I$,
$x(t)$ is governed by a controlled ordinary differential equation, 
and $i(t)$ is a Markov process with controlled transition probability matrix. 
See \cite{BenLi} for this application.

In this paper, we adopt the notion of viscosity solution to \erf{Pl}, for which
the reader may consult \cites{BaCa,BarB, CIL, CL, PLL}.

Now, we give our main assumptions on the system \erf{Pl}.
Throughout we implicitly assume that the functions $H_i$ are continuous in
$\T^n\tim\R^n\tim\R^m$.

We assume that $H$ is coercive, that is, for any $i\in\I$ and $R>0$,
\beq\tag{H1}\label{coercive}
\lim_{|p|\to\infty} \inf_{(x,u)\in\T^n\tim B_R^m}H_i(x,p,u) =\infty,
\eeq
where $B_R^m$ denotes the $m$-dimensional open ball with center at the origin and
radius $R$.

This is a standard assumption, under which any upper semicontinuous subsolution of  \erf{Pl} is Lipschitz continuous on $\T^n$.

We next assume that $H$ is convex in the variables $(p,u)$, that is,
\beq\tag{H2}\label{convex}
\text{for any $(x,i)\in\T^n\tim \I$, the function $(p,u)\mapsto H_i(x,p,u)$ is convex on $\R^n\tim\R^m$.}
\eeq

We assume that the Hamiltonian $H$ is monotone in the variable $u$,
that is, it satisfies
\beq\tag{H3}\label{monotone}
\left\{\
\begin{minipage}{0.8\textwidth} for any $(x,p)\in\T^n\tim\R^n$ and
$u=(u_i)_{i\in\I},\, v=(v_i)_{i\in\I}\in\R^m$,  if \ $u_k-v_k=\max_{i\in\I}(u_i-v_i)\geq 0$, \ then \
$H_k(x,p,u)\geq H_k(x,p,v)$.
\end{minipage}\right.
\eeq
This is a natural assumption implying that \erf{Pl}
should possess the comparison principle
between a subsolution and a supersolution.

When the coupling is linear, that is, when $H$ has the form
\beq\label{linear}
H_i(x,p,u)=G_i(x,p)+\sum_{j\in\I} b_{ij}(x)u_j\ \ \FOR i\in\I,
\eeq
condition \erf{monotone} is valid if and only if for each $x\in\T^n$, the matrix $B(x)$
is monotone matrix in the following sense
\beq \label{mono-matrix1}
b_{ij}(x)\leq 0 \ \ \IF i\not= j \ \ \ \AND \ \ \ \sum_{j\in\I}b_{ij}(x)\geq 0 \ \ \FORALL i\in\I,
\eeq
which is equivalent to that for each $x\in\T^n$, 
\beq \label{mono-matrix2}
\begin{minipage}{0.8\textwidth}
if $u=(u_i)\in\R^m$, $k\in\I$, and  
$u_k=\max_{i\in\I}u_i\geq 0$,\ then $(B(x)u)_k\geq 0$, 
\end{minipage}
\eeq
where $(B(x)u)_k$ denotes the $k$-th component of the $m$-vector $B(x)u$. We refer, for instance, to \cite[Lemma 3]{Icoupling} for the equivalence of the last two conditions \erf{mono-matrix1} and \erf{mono-matrix2}, while it is obvious 
in the linear coupling case \erf{linear} that 
\erf{monotone} and \erf{mono-matrix2} are equivalent each other.

When we deal with problem (P$_0$), we use the assumption that
\beq \tag{H4}\label{ergodic}
\begin{minipage}{0.8\textwidth} problem (P$_0$) has a solution in $C(\T^n)^m$.
\end{minipage}
\eeq
In the scalar case, that is, the case when $m=1$ and the case when $H(x,p,u)$ is independent of $u$, a natural problem which replaces 
(P$_0$) is the so-called ergodic problem that seeks a pair of a constant $c\in\R$ 
and a function $u\in C(\T^n)$ such that $u$ is a solution of 
\beq\label{ergo-gen}
H(x,Du)=c \ \ \IN \T^n.
\eeq
This problem is well-posed under \erf{coercive}, which means that there exists 
such a pair $(c,u)\in\R\tim \Lip(\T^n)$ and the constant $c$ is unique. 
See for this \cite{LPV}.  For such $(c,u)$, if we set $H_c=H-c$, then $u$ is 
a solution of (P$_0$), with $H$ replaced by $H_c$. If $v^\gl$ is a solution of 
\erf{Pl}, with current scalar Hamiltonian $H(x,p)$, then the function 
$v^\gl+\gl^{-1}c$ is a solution of \erf{Pl}, with $H$ replaced by $H_c$. 
This way, if $m=1$ and $H(x,p,u)$ is independent of $u$, 
then the vanishing discount problem can be transferred to the case where the limit 
problem (P$_0$) admits a solution.  Even in the scalar case and much more in the case where $m>1$, if $H(x,p,u)$ depends genuinely on $u$, then the reduction argument above does not work and the solutions of $H[u]=c$, 
where $c\not=0$, 
do not help investigate the vanishing discount problem for \erf{Pl}.

The rest of this paper is organized as follows.  In Section 2, we present
some basic properties concerning \erf{monotone} and a standard comparison and existence result of solutions of \erf{Pl} for $\gl>0$.
In Section 3, under additional hypotheses on the continuity of the Lagrangian $L$
and the compactness of the domain of $L$,
we study Green-Poisson measures for our system, which are crucial
in our asymptotic analysis. Section 4 establishes the compactness of the
support of the Green-Poisson measures and gives a representation theorem
for the solution of \erf{Pl}, with $\gl>0$, using the Green-Poisson measures.
We establish the main result for the vanishing discount
problem in Section 5. In Section 6, we establish an existence theorem for  the ergodic problem.

\section{Preliminaries}

We use the symbol $u\leq v$ (resp., $u\geq v$) for $m$-vectors $u,v\in\R^m$ to indicate
$u_i\leq v_i$ (resp., $u_i\geq v_i$) for all $i\in\I$.  Let $e_i$ denote the unit vector in $\R^m$ with unity as its $i$-th entry and let $\1$ 
denote the $m$-vector $(1,\ldots,1)\in\R^m$.

Concerning the monotonicity of $H$, we give a basic lemma.

\begin{lemma} \label{lem1}Assume that $H$ satisfies  \erf{monotone}. Let $\ga\geq 0$. 

\noindent 
\emph{(i)} For all  $i\in \I$ and $(x,p,u)\in\T^n\tim\R^n\tim\R^m$, 
\[
H_i(x,p,u+\ga\1)\geq H_i(x,p,u) .
\]
\emph{(ii)} For all  $i,j\in\I$ and $(x,p,u)\in\T^n\tim\R^n\tim\R^m$, if $i\not=j$, then 
\[
H_j(x,p,u+\ga e_i)\leq H_j(x,p,u)  
\]
\end{lemma}

\bproof (i) \ Fix $k\in\I$ and $(x,p,u)\in\T^n\tim\R^n\tim \R^m$. Set $v=u+\ga\1$ and note that
\[
(v-u)_k=\ga=\max_{i\in\I}(v-u)_i>0.
\]
By the monotonicity, we have
\[
H_k(x,p,v)\geq H_k(x,p,u).
\]
That is,
\[
H_k(x,p,u+\ga\1)\geq H_k(x,p,u). \qedhere
\]

(ii) \ Fix $i,j\in\I$ so that $i\not=j$. Let $(x,p,u)\in\T^n\tim\R^n\tim\R^m$.
Note that $u_k-(u+\ga e_i)_k=0$ if $k\not=i$ and $=-\ga<0$ if $k=i$, which can be stated that $u_j-(u+\ga e_i)_j=\max_{k\in\I}[u_k-(u+\ga e_i)_k]=0$. By \erf{monotone}, we have
\[
H_j(x,p,u)\geq H_j(x,p,u+\ga e_i). \qedhere
\]
\eproof

The following theorem is well-known: see \cite{EL, IK}
 for instance for a general background and \cite[the proof of Theorem~1]{Icoupling}
for some details how to adapt general results to \erf{Pl}.

\begin{theorem}\label{thm2-1} Assume \eqref{coercive} and \eqref{monotone}. Let $\gl>0$. Then there exists a unique
solution $v^\gl\in\Lip(\T^n)^m$ of \erf{Pl}. Also, if $v=(v_i), w=(w_i)$ are,
respectively, upper and lower semicontinuous on $\T^n$ and a subsolution and
a supersolution of \erf{Pl}, then $v\leq w$ on $\T^n$.
\end{theorem}

We remark that if $H$ satisfies \erf{coercive} and \erf{monotone}, then so does the 
Hamiltonian $H_\ep(x,p,u):=H(x,p,u)+\ep u$, with $\ep>0$ and that if  our 
problem \erf{Pl} has $H_\ep$ in place of $H$, then the 
limit problem (P$_0$) reads $\ep u +H[u]=0$ in $\T^n$ and has a unique solution 
due to Theorem \ref{thm2-1}. Thus, the asymptotic analysis for the vanishing 
discount problem in such cases is fairly easy.

With reference to \cite{IK}, we outline the proof of the theorem above.

\bproof[Outline of proof] We choose a constant $C>0$ so that
\[
\max_{(x,i)\in\T^n\tim \I}|H_i(x,0,0)|\leq C.
\]
Note by Lemma \ref{lem1} that
\[
H_i(x,0,-\gl^{-1}C\1)\leq H_i(x,0,0)\leq H_i(x,0,\gl^{-1}C\1) \ \ \FOR x\in\T^n.
\]
By using this, it is easily checked that the functions
$f(x)=\gl^{-1}C\1$ and $g(x)=-\gl^{-1}C\1$ are
a supersolution and subsolution of \erf{Pl} and satisfy $f\geq g$ on $\T^n$.

Our assumption \erf{monotone} implies the quasi-monotonicity of $H$ in \cite{IK}
as was shown in \cite[Lemma 4.8]{IK}. By \cite[Theorem3.3]{IK}, the function
$z=(z_i)_{i\in\I}$ on $\T^n$ given by
\[
z_i(x)=\sup\{\zeta_i(x)\mid g\leq \zeta \leq f \IN \T^n,\
\zeta \text{ is a subsolution of }\erf{Pl}\} \ \ \FOR (x,i)\in\T^n\tim\I,
\]
is a solution of \erf{Pl}, in the sense that $z^*=(z^*_i)_{i\in\I}$ and $z_*=(z_{i*})_{i\in\I}$,  where each $z_i^*$ and $z_{i*}$ are respectively
the upper and lower semicontinuous envelope of $z_i$, are respectively a subsolution and a supersolution of \erf{Pl}.

By the definition of $z$, it is easy to infer that for all $i\in\I$, $z_i=z_i^*$ and the function $z_i$ is upper semicontinuous on $\T^n$. By \erf{coercive} (the coercivity of $H$),
we deduce that the function $z$ is Lipschitz
continuous on $\T^n$.

To see that the comparison between $v$ and $w$,
we apply \cite[Theorem 4.7]{IK} to $v$ and $z$ as well as $z$ and $w$, to conclude that
$v\leq z$ and $z\leq w$ on $\T^n$, which implies that $v\leq w$ on $\T^n$. Here, the comparison theorem
 (\cite[Theorem 4.7]{IK})
requires the regularity of $H$ (see \cite[(A.2)]{IK}),
which can be reduced just to the continuity of $H$
since $z$ is Lipschitz continuous on $\T^n$.
This reduction of regularity of $H$ is a standard observation and we leave it to the interested reader
to adapt the proof of  \cite[Theorem 4.7]{IK} to this
case.
\eproof

Setting
\[
L_i(x,\xi,\eta):=\sup_{(p,u)\in\R^n\tim\R^m}[\xi\cdot p+\eta\cdot u-H_i(x,p,u)] \ \ \FOR (x,i,\xi,\eta)\in\T^n\tim\I\tim\R^n\tim\R^m,
\]
by the convex duality we have
\[
H_i(x,p,u)=\sup_{(\xi,\eta)\in\R^n\tim\R^m}[\xi\cdot p+\eta\cdot u-L_i(x,\xi,\eta)] \ \ \FOR (x,i,p,u)\in\T^n\tim\I\tim\R^n\tim\R^m.
\]
We call $L_i$ (resp., $(L_i)_{i\in\I}$) the Lagrangian of $H_i$ (resp.,  the Lagrangian of $(H_i)_{i\in\I}$).
Similarly, we call $H_i$ (resp.,  $(H_i)_{i\in\I}$) the Hamiltonian of $L_i$ (resp.,  the Hamiltonian of $(L_i)_{i\in\I}$).  

It should be remarked that the functions $L_i$ are 
lower semicontinuous on $\T^n\tim\R^n\tim\R^m$.   

We examine the Lagrangian in the linear coupling case \erf{linear}. 
By the definition of $L_i$ and a simple manipulation, we deduce that
\beq\label{linear-Lag} \bald
L_i(x,p,\xi,\eta)
&=\sup_{p\in\R^n}(\xi\cdot p-G_i(x,p))+\sup_{u\in\R^m}(\eta\cdot u-(B(x)u)_i)
\\&=\sup_{p\in\R^n}(\xi\cdot p-G_i(x,p))+\0_{\{b_i(x)\}}(\eta),
\eald\eeq
where for $x\in\T^n$, $b_i(x):=(b_{ij})_{j\in\I}\in\R^m$ and, for any sets $X\subset Y$, 
$\0_X$ denotes the indicator function of $X$ on $Y$ given by 
\[
\0_X(y)=
\bcases
0   & \IF\  y\in X, \\[3pt]
+\infty\ & \IF \ y\in Y\stm X.
\ecases
\]

\begin{lemma}\label{lowerB} Assume \erf{coercive}--\erf{convex}. 

\noindent 
\emph{(i)}~We have
\beq\label{lB1}
L_i(x,\xi,\eta)\geq -H_i(x,0,0) \ \ \FOR  (x,i,\xi,\eta)\in \T^n\tim \I\tim\R^n\tim \R^m.
\eeq
\emph{(ii)}~For any $A>0$ there exists a constant $C_A$ such that
\beq
L_i(x,\xi,\eta)\geq A|\xi|-C_A \ \ \FOR (x,i,\xi,\eta)\in \T^n\tim\I\tim \R^n\tim \R^m.
\eeq
\end{lemma}

\bproof Fix $i\in\I$. We have
\[
L_i(x,\xi,\eta)=\sup_{(p,u)\in\R^n\tim\R^m}(\xi\cdot p+\eta\cdot u-H_i(x,p,u))
\geq -H_i(x,0,0),
\]
and
\[
L_i(x,\xi,\eta)=\sup_{(p,u)\in\R^n\tim\R^m}(\xi\cdot p+\eta\cdot u-H_i(x,p,u))
\geq A|\xi|-H_i(x,A\xi/|\xi|,0) \ \ \IF \xi\not=0.
\]
Hence, setting
\[
C_A=\sup_{(x,i,\xi)\in\T^n\tim\I\tim B_A^n}H_i(x,\xi,0),
\]
we obtain
\[
L_i(x,\xi,\eta)\geq A|\xi|-C_A. \qedhere
\]
\eproof

Lemma \ref{lowerB}, (ii) asserts that the functions $L_i(x,\xi,\eta)$ have
a superlinear growth as $|\xi| \to \infty$.

We give a characterization of the monotonicity \erf{monotone} of $(H_i)_{i\in\I}$
through $(L_i)_{i\in\I}$. For $k\in\I$, we write
\[
Y_k:=\{(\eta_i)_{i\in\I}\in\R^m\mid \eta_i\leq 0 \ \IF i\not=k,\, \sum_{i\in\I}\eta_i\geq 0\},
\]
and \ $
\dom L_k=\{(x,\xi,\eta)\mid L_k(x,\xi,\eta)<\infty\}.$

\begin{\prop} \label{prop2.1} Assume \erf{coercive}--\erf{convex}. Then $\Hall$ satisfies
\erf{monotone} if and only if
\[
\dom L_i \subset \T^n\tim\R^n\tim Y_i \ \ \FORALL i\in\I.
\]
\end{\prop}

One can check directly that, in the linear coupling case \erf{linear}, 
if \erf{coercive}--\erf{monotone} hold, then the inclusion above is valid. Indeed, in this case, the coupling matrix $B(x)=(b_{ij}(x))$ satisfies \erf{mono-matrix1}, which 
implies that $(b_{ij}(x))_{j\in\I}\in Y_i$ for all $(x,i)\in\T^n\tim\I$, 
and, by \erf{linear-Lag}, we conclude that $\dom L_i\subset \T^n\tim\R^n\tim Y_i $ for all $i\in\I$.

\bproof  We assume first that $\Hall$ satisfies \erf{monotone}.
Fix any $(x,k,\xi,\eta)\in\T^n\tim\I\tim\R^n\tim\R^m$ and suppose that $\eta=(\eta_i)_{i\in\I}\not\in Y_k$. We have either $\eta_j>0$ for some $j\not=k$ or
$\sum_{i\in\I}\eta_i<0$.

Consider the case when $\eta_j>0$ for some $j\not=k$. Let $t>0$,  and,  
by (ii) of Lemma \ref{lem1}, 
we have
\[
H_k(x,p,0)\geq H_k(x,p,te_j) \ \ \FOR p\in\R^n,
\]
and hence,
\[
\xi\cdot p+\eta\cdot te_j -H_k(x,p,te_j)
\geq \xi\cdot p+t\eta_j -H_k(x,p,0) \ \ \FOR p\in\R^n,
\]
which implies that $L_k(x,\xi,\eta)=\infty.$

Consider next the case when $\sum_{i\in\I}\eta_i<0$. For $t>0$, we 
observe 
by (i) of Lemma \ref{lem1} that
$H_k(x,p,0)\geq H_k(x,p,t\1)$ for all $p\in\R^n$. Consequently,
\[
\xi\cdot p+\eta\cdot t\1 -H_k(x,p,t\1)
\geq \xi\cdot p-t\sum_{i\in\I}\eta_i -H_k(x,p,0) \ \ \FOR p\in\R^n,
\]
which shows that $L_k(x,\xi,\eta)=\infty$.  We thus conclude that
$\dom L_k\subset \T^n\tim\R^n\tim Y_k$.

Next, we assume that $\dom L_i\subset \T^n\tim\R^n\tim Y_i$ for all $i\in\I$.
It is obvious that for any $(x,i,p,u)\in\T^n\tim\I\tim\R^n\tim\R^m$,
\[
H_i(x,p,u)=\sup_{(\xi,\eta)\in \R^n\tim Y_i}[\xi\cdot p+\eta\cdot u-L_i(x,\xi,\eta)].
\]

Fix any $(x,p)\in\T^n\tim\R^n$ and $u,v\in\R^m$. Assume that for some $k\in\I$,
\[
(u-v)_k=\max_{i\in\I}(u-v)_i\geq 0,
\]
which can be stated as
\[
(u-v)_k-(u-v)_i\geq 0 \ \ \FOR i\in\I \ \ \AND  \ \ (u-v)_k\geq 0.
\]
Let $\eta=(\eta_i)_{i\in\I}\in Y_k$. Multiplying the first inequality above by $\eta_i$, with $i\not=k$, we get
\[\bald
0&\,\geq \sum_{i\not=k}\eta_i[(u-v)_k-(u-v)_i]
=\sum_{i\in\I}\eta_i[(u-v)_k-(u-v)_i]
\\&\,=(u-v)_k\sum_{i\in\I}\eta_i-\sum_{i\in\I}\eta_i(u-v)_i.
\eald\]
Since $(u-v)_k\geq 0$ and $\sum_{i\in\I}\eta_i\geq 0$, we infer from the above that
$\eta\cdot u\geq \eta\cdot v$.  Thus, we have
\[\bald
H_k(x,p,u)&\,=\sup_{(\xi,\eta)\in\R^n\tim Y_k}
[\xi\cdot p+\eta\cdot u-L_k(x,\xi,\eta)]
\\&\,\geq
\sup_{(\xi,\eta)\in\R^n\tim Y_k}
[\xi\cdot p+\eta\cdot v-L_k(x,\xi,\eta)]=H_k(x,p,v).
\eald
\]
This shows that $\Hall$ is monotone, which completes the proof.
\eproof

\section{Green-Poisson measures: in a regular case}

In what follows,  given a topological space $X$, $\cB(X)$ denotes the $\gs$-algebra of Borel sets in $X$, and $\M(X)$ and $\M_+(X)$ denote, respectively,  the spaces 
of Borel measures on $X$ having bounded variation and of nonnegative finite Borel
measures on $X$. Also, $C_\rmb(X)$ denotes the space of 
bounded continuous functions on $X$. 

For any $\nu\in\M(\T^n\tim \R^n\tim\R^m)$ and
integrable function $\phi$ on $\T^n\tim \R^n\tim\R^m$ with respect to $\nu$, we write
\[
\du{\nu,\phi}=\int_{\T^n\tim \R^n\tim\R^m}\phi(x,\xi,\eta)\,\nu(dx d\xi d\eta).
\]
Similarly, for any $\nu=(\nu_i)_{i\in\I}\in \M(\T^n\tim\R^n\tim\R^m)^m$ and Borel function
$\phi=(\phi_i)_{i\in\I}$ on $\T^n\tim \R^n\tim\R^m$,
we write
\[
\du{\nu,\phi}=\sum_{i\in\I}\du{\nu_i,\phi_i}
\]
if $\phi_i$ is integrable with respect to $\nu_i$ for every $i\in\I$.

For $\gl>0$, we define the function $S^\gl : \R^m\to\R$ by  
$S^\gl(\eta)=\gl+\sum_{i\in\I}\eta_i$ for $\eta=(\eta_i)_{i\in\I}\in \R^m$, 
and,  we write $\P^\gl$ for the set of all $\mu=(\mu_i)_{i\in\I}\in \M_+(\T^n\tim \R^n\tim\R^m)^m$ such that
\beq\label{normal>0}
\du{\mu_i,|\xi|+|\eta|}<\infty 
\ \ \FORALL i\in \I \ \ \ \AND \ \ \ \du{\mu,S^\gl\,\1}=1,
\eeq
where $|\xi|+|\eta|$ denotes the function: $\R^n\tim\R^m\ni (\xi,\eta) \mapsto |\xi|+|\eta|\in\R$. Note that $S^\gl\1$ is the function: $\R^m\ni\eta \mapsto S^\gl(\eta)\1\in\R^m$, which can be regarded as a function of $(x,\xi,\eta)\in\T^n\tim\R^{n+m}$. 
We write $\P^0$ for the set of all $\mu=(\mu_i)_{i\in\I}\in \M_+(\T^n\tim \R^n\tim\R^m)^m$ such that
\beq\label{normal=0}
\du{\mu_i,|\xi|+|\eta|}<\infty \ \ \FORALL i\in\I \ \ \ \AND \ \ \ 
\du{\mu,\,\1}\leq 1. 
\eeq

Similarly to the standard definition \cite{DFIZ, DZ2} (see also \cite{Man,Mat}) of Mather measures, we introduce the closed measures as follows. 
We call any $\mu\in\P^0$ a closed measure (associated with $\gl=0$) provided it satisfies 
\beq\label{gpr2}
\du{\mu, \xi\cdot D\psi+\eta\cdot\psi\, \1}=0  \ \ \text{ for all }\psi\in C^1(\T^n)^m\quad\IF \ \gl=0,
\eeq
and denote by $\frak{C}(0)$ the set of all such closed measures associated with $\gl=0$. 
In \erf{gpr2} above, and henceforth, we use the notation that the functions: 
$\T^n\tim \R^n\ni (x,\xi)\mapsto (\xi\cdot D\psi_i(x))_{i\in\I}\in\R^m$ and $\T^n\tim \R^m\ni (x,\eta) \mapsto \eta\cdot \psi(x)\in\R$ are denoted 
by $\xi\cdot D\psi$ and $\eta\cdot \psi$, respectively.
For $(z,k,\gl)\in\T^n\tim\I\tim(0,\,\infty)$, we call any $\mu\in\P^\gl$ a closed measure (associated with $(z,k,\gl)$) provided it satisfies  
\beq \label{gpr1}
\du{\mu,\xi\cdot D\psi+\eta\cdot \psi\,\1+\gl\psi}=\psi_k(z) \ \ \text{ for all }
\psi=(\psi_i)_{i\in\I}\in C^1(\T^n)^m \quad\IF \ \gl>0.
\eeq 
The set of all closed measures associated with $(z,k,\gl)$ is denoted by $\frak{C}(z,k,\gl)$.

We now introduce the following working hypothesis.
\beq\tag{H5}\label{compact-continuous} \left\{
\begin{minipage}{0.9\textwidth}There exist nonempty, compact, convex sets $K_1\subset \R^n$ and $K_2\subset \R^m$ such that 
$K_1$ is a neighborhood of the origin and such that 
for $i\in\I$,
\[\bald
L_i &\,\in C(\T^n\tim K_1\tim(K_2\cap Y_i)), \\
H_i(x,p,u)&\,=\sup_{(\xi,\eta)\in K_1\tim(K_2\cap Y_i)}
[\xi\cdot p+\eta\cdot u-L_i(x,\xi,\eta)]\ \ \FOR (x,p,u)\in\T^n\tim\R^n\tim\R^m.
\eald\]
\end{minipage}\right.
\eeq

\def\0{{\mathbf 0}}

Assuming \erf{compact-continuous} in addition, we have
\[
L_i(x,\xi,\eta)=L_i(x,\xi,\eta)+
\0_{K_1\tim (K_2\cap Y_i)}(\xi,\eta)
\ \ \FOR (x,\xi,\eta)\in\T^n\tim\R^n\tim\R^m.
\]

We remark that under \erf{compact-continuous}, the functions
$H_i(x,p,u)$ grows at most  linearly as $|(p,u)| \to \infty$.

\begin{theorem}\label{gpr-thm2} Assume \erf{coercive}--
\erf{monotone} and \erf{compact-continuous}. Let $(z,k,\gl) \in \T^n\tim\I\tim
(0,\,\infty)$ and let $v^\gl=(v^\gl_i)_{i\in\I}\in C(\T^n)^m$ be the solution of \erf{Pl}.
Then there
exists $\mu\in\frak{C}(z,k,\gl)$ such that
\beq\label{gpr3}
v_k^\gl(z)=\du{\mu,L}=\min_{\nu\in\frak{C}(z,k,\gl)}\du{\nu,L}.
\eeq
\end{theorem}

Remark that, thanks to Proposition \ref{prop2.1},  if $\Hall$ satisfies
\erf{coercive}, \erf{convex} and \erf{compact-continuous},
then it has the property \erf{monotone} as well.

The theorem above and Theorem \ref{gpg-thm1}, stated later, are generalizations 
of the previous results in \cite{Go08, CaCDGo, DZ2, DFIZ, MT, GMT, IMT1, IMT2, Icoupling}

Our proof of Theorem \ref{gpr-thm2} is close to the ones in \cite{IMT1,IMT2} 
in the technicality and is similar to the ones in \cite{Go, IS} in the use of duality. It 
depends crucially on a minimax theorem, and, in the application of the minimax theorem, it is essential to make the set compact 
on which the measures $\mu\in \frak{C}(z,k,\gl)$, having the property $\du{\mu, L}<\infty$, are supported. 
The condition \erf{compact-continuous} realizes all such measures $\mu$ to be supported  on the compact set $\prod_{i\in\I}\T^n\tim K_1\tim (K_2\cap Y_i)$, which 
we mean that $\supp \mu_i\subset \T^n\tim K_1\tim (K_2\cap Y_i)$ for all $i\in\I$.

In the next section, we remove the restriction \erf{compact-continuous} on $L$ 
adopted in Theorem \ref{gpr-thm2} by appealing the fact that the solution $v^\gl$ of \erf{Pl} 
is Lipschitz continuous.

We call a \emph{Green-Poisson measure}
any measure $\mu\in\frak{C}(z,k,\gl)$ that is a minimizer of the most right hand of \erf{gpr3}.

\def\tilF{\cF}

\def\tilG{\cG}

For the proof of Theorem \ref{gpr-thm2}, we 
assume henceforth \erf{coercive}--\erf{monotone} and \erf{compact-continuous}, and introduce the following notation. 
Set $Z_i=K_1\tim(K_2\cap Y_i)$ for $i\in\I$ and note that for any $i\in\I$, $Z_i$ is a compact convex subset of $\R^n\tim\R^m$. 
Let $\gl>0$ and $\tilF(\gl)$ denote the set of all $(\phi,u)\in \prod_{i\in\I}C(\T^n\tim Z_i)\tim C(\T^n)^m$ such that for any $(x,i)\in\T^n\tim \I$, $\phi_i(x,\xi,\eta)$ is convex 
in the variable $(\xi,\eta)$ on $Z_i$ and 
such that $u$ is a subsolution of $\gl u+H_\phi[u]=0$ in $\T^n$,
where $H_\phi=(H_{\phi,i})_{i\in\I}$ is given by
\beq\label{def_phi_i}
H_{\phi,i}(x,p,u)=\max_{(\xi,\eta)\in Z_i}(\xi\cdot p+\eta\cdot u-\phi_i(x,\xi,\eta)).
\eeq

Here, we note by the compactness of $Z_i$ that $H_{\phi,i}$ is continuous on 
$\T^n\tim\R^n\tim\R^m$. 
Also, if we identify $\phi_i$ with the function $\tilde\phi_i$ on $\T^n\tim\R^n\tim\R^m$ given by
\beq \label{ext_to_whole}
\tilde\phi_i(x,\xi,\eta)=
\bcases
\phi_i(x,\xi,\eta) \ &\IF \ (\xi,\eta)\in Z_i,
\\
+\infty \ &\text{ otherwise}, 
\ecases
\eeq
then \erf{def_phi_i} reads 
\[
H_{\phi,i}(x,p,u)=\max_{(\xi,\eta)\in \R^n\tim\R^m}(\xi\cdot p+\eta\cdot u-\phi_i(x,\xi,\eta)).
\]
Notice that the functions $\tilde\phi_i(x,\xi,\eta)$
defined by \erf{ext_to_whole} are lower semicontinuous in the variable $(x,\xi,\eta)$ and convex in the variable $(\xi,\eta)$.  By the convex duality, we see that $\tilde\phi_i$ is the Lagrangian 
of $H_{\phi,i}$.  Arguing similarly to the proof of (ii) of Lemma \ref{lowerB}, with 
$(\tilde\phi_i, H_{\phi,i})$ in place of $(H_i,L_i)$, we see that $H_\phi$ satisfies 
\erf{coercive}. 
By Proposition \ref{prop2.1}, we easily see that $H_\phi$ satisfies \erf{monotone}. 
Moreover, it is clear that $H_\phi$ satisfies \erf{compact-continuous}.

In what follows, for any $(\phi,u)\in \cF(\gl)$, we identify the function $\phi_i$ on $\T^n\tim Z_i$ with $\tilde\phi_i$ according to the situation.  In particular, if $v^\gl\in C(\T^n)$ is the solution of \erf{Pl}, then we have \ $(L,v^\gl)\in\cF(\gl)$. 
Note in addition that if $\phi=0\in \prod_{i\in\I}C(\T^n\tim Z_i)$, then 
$\phi_i(x,\xi,\eta)$ is convex in $(\xi,\eta)$ for all $i\in\I$ and $H_{\phi}(\cdot,0,0)=0$.  
Hence, $(0,0)\in\cF(\gl)$.

For $(z,k,\gl)\in\T^n\tim\I\tim(0,\,\infty)$, we set
\[\bald
\tilG(z,k,\gl)&\,=\{\phi-u_k(z)S^\gl(\eta)\1\mid (\phi,u)\in\tilF(\gl)\},\\
\P^\gl(K_1,K_2)&\,= \{\mu\in\P^\gl\mid 
\supp \mu_i\subset \T^n\tim Z_i \ \FORALL i\in\I\}, 
\\
\tilG\pr(z,k,\gl)&\,=\{\mu\in\P^\gl(K_1,K_2)\mid 
\du{\mu,f}\geq 0 \ \FORALL f\in\tilG(z,k,\gl)\}.
\eald
\]

We recall that, by definition, the support of measure $\mu\in \M(X)$ for a topological space $X$ is defined as the closed set
\[
\supp \mu =X\stm \bigcup\{U\subset X \mid |\mu|(U)=0,\, U\text{ is open}\},
\]
where $|\mu|$ denotes the total variation of the measure $\mu$.  
Accordingly, if $X$ has a countable basis of its topology, we have
\[
\mu(X\stm \supp \mu)=0, 
\]
and 
\[
\int_X \phi(x)\mu(dx)=\int_{\supp \mu} \phi(x)\mu(dx) \ \ \FORALL \phi\in C_\rmb(X),\,\mu\in\M(X). 
\]
In what follows, when $Q\subset X$ is a closed subset of a topological space $X$ with a countable basis, we may identify $\mu\in\M(X)$ satisfying 
$\supp \mu\subset Q$ with its restriction $\mu|_Q$ to $Q$, defined by 
\[
\mu|_Q(A)=\mu(A) \ \ \FORALL A\in\cB(Q). 
\] 
Then we have
\beq\label{restriction}
\int_X \phi(x)\mu(dx)=\int_{Q} \phi(x)\,\mu|_Q(dx) \ \ \FORALL \phi\in C_\rmb(X). 
\eeq

Recalling that $(L,v^\gl)\in \cF(\gl)$, where $v^\gl$ is a solution of \erf{Pl},  
we easily infer that
\[
v^\gl_k(z)\leq \du{\mu,L} \ \ \FORALL \mu\in \cG\pr(z,k,\gl). 
\]

Since the functions $L_i$ 
are bounded from below, for any $\mu\in\P^\gl$, the inequality 
$\du{\mu,L}<\infty$ always makes sense and, by \erf{compact-continuous}, 
we have $\du{\mu,L}<\infty$ if and only if $\supp\mu_i\subset \T^n\tim Z_i$ for all $i\in\I$.  Hence, we have 
\beq \label{Two_chara_Pl}
\P^\gl(K_1,K_2)=\{\mu\in\P^\gl\mid \du{\mu,L}<\infty\}. 
\eeq

\begin{lemma} \label{cone} The set $\tilF(\gl)$ is a convex cone in
$\prod_{i\in\I}C(\T^n\tim Z_i)\tim C(\T^n)^m$ with vertex
at the origin.
\end{lemma}

\bproof Recall \cite[Remark 2.5]{Bar} that for any $u\in \Lip(\T^n)^m$,
$u$ is a subsolution of
\[
\gl u+H[u]=0 \ \ \IN\T^n
\]
if and only if for all $i\in\I$,
\[
\gl u_i(x)+H_i(x,Du_i(x),u(x))\leq 0 \ \ \text{ a.e. in }\T^n,
\]
and by the coercivity \erf{coercive} that for any $(\phi,u)\in\tilF(\gl)$, we have
$u\in\Lip(\T^n)^m$.

Fix $(\phi,u), (\psi,v)\in\tilF(\gl)$ and $t,s\in[0,\infty)$. Fix $i\in\I$ and
observe that
 \[\bald
&\,\gl u_i(x)+H_{\phi,i}(x,Du_i(x),u(x))\leq 0 \ \ \text{ a.e. in }\T^n,
\\&\,\gl v_i(x)+H_{\psi,i}(x,Dv_i(x),v(x))\leq 0 \ \ \text{ a.e. in }\T^n,
\eald
\]
which imply that there is a set $N\subset\T^n$ of Lebesgue measure zero such that
\[\bald
&\,\gl u_i(x)+\xi\cdot Du_i(x) +\eta\cdot u(x)\leq \phi_i(x,\xi,\eta)  \ \ \text{ for all } (x,\xi,\eta)
\in(\T^n\stm N)\,\tim\,Z_i,
\\&\,\gl v_i(x)+\xi\cdot Dv_i(x)+\eta\cdot u(x)\leq \psi_i(x,\xi,\eta) \ \ \text{ for all }(x,\xi,\eta)\in(\T^n\stm N)\,\tim\,Z_i.
\eald
\]
Multiplying the first and second by $t$ and $s$, respectively, adding
the resulting inequalities and setting $w=tu+sv$, we obtain
\[
\gl w_i(x)+ \xi\cdot Dw_i(x)+\eta \cdot w(x) \leq
(t\phi_i+s\psi_i)(x,\xi,\eta)  \ \ \text{ for all } (x,\xi,\eta)
\in(\T^n\stm N)\,\tim\,Z_i,
\]
which implies that $t(\phi,u)+s(\psi,v)\in\tilF(\gl)$.
\eproof

\begin{lemma} \label{gpr-thm1}  
Let $(z,k,\gl)\in\T^n\tim\I\tim(0,\,\infty)$ and $\mu=(\mu_i)_{i\in\I}\in \P^\gl(K_1,K_2)$. Then, we have 
$\mu\in \frak{C}(z,k,\gl)$ if and only if $\mu\in\cG\pr(z,k,\gl)$.  
\end{lemma}

\bproof  Assume first that $\mu\in\frak{C}(z,k,\gl)$. 
Fix any $(\phi,u)\in\cF(\gl)$ 
and recall that, in the viscosity sense,
\[ 
\gl u+H_\phi[u]\leq 0 \ \ \IN \T^n.  
\]
Thanks to the coercivity property \erf{coercive} of $H_\phi$, $u$ is Lipschitz continuous on $\T^n$.  
In view of the continuity of $H_\phi$ and the convex property \erf{convex} 
of $H_\phi$, mollifying $u$, we may choose, for each $\ep>0$,  
a function $u^\ep\in C^1(\T^n)^m$ such that 
$\gl u^\ep+H_\phi[u^\ep]\leq \ep\,\1$ in $\T^n$ and $\|u-u^\ep\|_\infty<\ep$.  Hence, we have
\[
\gl u^\ep+\xi\cdot Du^\ep(x)+\eta\cdot u^\ep(x)\,\1 \leq \phi(x,\xi,\eta)+\ep\,\1.
\]
Using \erf{gpr1} and integrating the inequality above with respect to $\mu$, we obtain
\[\bald
u^\ep_k(z)&\,=\du{\mu,\, \xi\cdot Du^\ep+\eta\cdot u^\ep\,\1+\gl u^\ep}
\leq \du{\mu,\, \phi+\ep\,\1}
\eald
\]
and, after sending $\ep\to 0$,
\[
0\leq \du{\mu,\phi}-u_k(z)=\du{\mu, \phi-u_k(z)S^\gl \1},
\]
which implies, together with the assumption that $\supp \mu_i\subset 
\T^n\tim Z_i$ for all $i\in\I$, that \ $\mu\in\cG\pr(z,k,\gl)$.

Next, we assume that  \ $\mu\in\cG\pr(z,k,\gl)$. Fix any $\psi \in C^1(\T^n)^m$, 
set $\phi=\xi\cdot D\psi+\eta\cdot \psi +\gl\psi$, which is a function on $\T^n\tim\R^{n+m}$, and 
observe that $\gl \psi+H_\phi[\psi]\leq 0$ in $\T^n$, i.e., $(\phi,\psi)\in\cF(\gl)$.  
Hence, by the definition of $\cG\pr(z,k,\gl)$, we have  
\[
0\leq \du{\mu, \phi-\psi_k(z)S^\gl\1}=\du{\mu,\xi\cdot D\psi+\eta\cdot \psi \1 
+\gl \psi}-\psi_k(z).
\]
The inequality above holds also for $-\psi$ in place of $\psi$, which reads
\[
0\geq \du{\mu,\xi\cdot D\psi+\eta\cdot \psi \1 
+\gl \psi}-\psi_k(z).
\]
Thus, we have 
\[ 
\du{\mu,\xi\cdot D\psi+\eta\cdot \psi \1 
+\gl \psi}=\psi_k(z),
\]
and conclude that $\mu\in\frak{C}(z,k,\gl)$. 
\eproof

\begin{lemma} \label{Dirac_measure} Let $i\in\I$, $\gl>0$, 
and $(\bar x,\bar\xi,\bar\eta)\in\T^n\tim Z_i$, 
and let $\gd_{(\bar x,\bar\xi,\bar\eta)}$ denote the Dirac measure at $(\bar x,\bar\xi,\bar\eta)$. Then, 
$\left(S^\gl\right)^{-1}\gd_{(\bar x,\bar\xi,\bar\eta)}e_i$ is 
a member of $\P^\gl(K_1,K_2)$. 
\end{lemma}

\bproof Note first that $S^\gl(\eta)\geq \gl>0$ for all $\eta\in K_2\cap Y_i$. 
It follows immediately that $\left(S^\gl\right)^{-1}\gd_{(\bar x,\bar\xi,\bar\eta)}\in\M^+(\T^n\tim Z_i)$, $\supp S^\gl\gd_{(\bar x,\bar\xi,\bar\eta)}=\{(\bar x,\bar \xi,\bar\eta)\}\subset \T^n\tim Z_i$, 
\[
\du{\left(S^\gl\right)^{-1}\gd_{(\bar x,\bar\xi,\bar\eta)}e_i,\,(|\xi|+|\eta|)\1}
=\left(S^\gl(\bar\eta)\right)^{-1}(|\bar \xi|+|\bar\eta|)<\infty,
\]
and \ $
\du{\left(S^\gl\right)^{-1}\gd_{(\bar x,\bar\xi,\bar\eta)}e_i,\,S^\gl\1}
=1.$ 
Thus, we we see that $\left(S^\gl\right)^{-1}\gd_{(\bar x,\bar\xi,\bar\eta)}e_i\in \P^\gl(K_1,K_2)$. 
\eproof 

For the reader's convenience, we state a minimax theorem (\cite[Corollary 2]{Te}).  

\begin{proposition} \label{minimax} Let $K$ and $Y$ be convex subsets of vector
spaces. Assume in addition that $K$ is a compact space.  Let $f\mid K \tim Y \to \R$ be
a function satisfying: 

\noindent 
\emph{(i)} For each $y \in Y$, the function: $x \mapsto f(x,y)$ is lower semicontinuous and
convex on $K$. 

\noindent 
\emph{(ii)} For each $x \in K$, the function: $y \mapsto f(x,y)$ is concave on $Y$.

\noindent Then
\[
\sup_{y\in Y}\min_{x\in K}f(x,y)=\min_{x\in K}\sup_{y\in Y}f(x,y).
\]
\end{proposition}  

We remark that in \cite[Corollary 2]{Te}, it is assumed that $K$ is a convex compact subset of a  topological vector space $X$, but in its proofs, the compatibility 
of the linear structure and the topological structure (i.e., the continuity 
of addition and scalar multiplication) of $X$ is not used and the proposition above is valid.

In the application below of Proposition \ref{minimax}, we take $K$ to be a bounded subset of the Banach space $\M(\gS)$ with the total variation norm, where 
$\gS$ is a compact subset of $\T^n\tim\R^{n+m}$.  

Let $\gS$ be a compact subset of $\T^n\tim\R^{n+m}$. By the Riesz representation theorem,  for each $F\in C(\gS)^*$, there exists a unique (regular) Borel measure 
$\mu$ on $\gS$ such that for all $\phi\in C(\gS)$, 
\[
F(\phi)=\int_{s\in \gS}\phi(s)\mu(ds)=\du{\mu,\phi}. 
\]
The mapping $\iota_\gS$ of $F\in C(\gS)^*$ to $\mu\in\M(\gS)$, given above, is an 
isomorphism between two Banach spaces.  Through the mapping $\iota_\gS: C(\gS)^* \to \M(\gS)$, 
the weak star convergence corresponds to the weak convergence of measures. 
 
Thanks to the Banach-Alaoglu theorem, we know that any closed ball $B$ (in the strong topology) of $C(\gS)^*$, equipped with the weak star topology, is a compact metrizable space.  Moreover, if $B$ is such a ball and $N$ is a closed subset (in the weak star topology) of $B$, then $N$ is a compact subset of $B$.  
These say that if $D$ is a closed ball (in the total variation norm) of 
$\M(\gS)$, then $D$ is a compact metrizable space with the topology of the weak convergence of measures and so is any $K\subset D$ that is sequentially closed in the topology of the weak convergence of measures.

The following lemma is a simple consequence of the discussion above. 

\begin{lemma} \label{PS} Let $a>0$ and $\gS=(\gS_i)_{i\in\I}$ be a 
collection of compact subsets $\gS_i$ of $\T^n\tim\R^{n+m}$.
Let $\P(\gS,a)$ denote the collection of $\mu=(\mu_i)\in 
\M_+(\T^n\tim\R^{n+m})^m$ such that $\supp \mu_i\subset \gS_i$ 
for all $i\in\I$ and such that $\du{\mu,\1}\leq a$.       
Then, $\P(\gS,a)$ is a compact metrizable space  
with the topology of weak convergence of measures. 
\end{lemma} 

It is to be noticed that in the lemma above, $\P(\gS,a)$ is sequentially compact.  

\bproof In view of \erf{restriction}, it is clear that any sequence of 
measures $\mu^q=(\mu_i^q)\in\P(\gS,a)$ on $\T^n\tim\R^{n+m}$, with $q\in\N$, 
converges to $\mu=(\mu_i)$ weakly in the sense of measures 
(i.e., in the topology of weak convergence of measures) if and only if, for each $i\in\I$,
the sequence of measures $\mu_i^q|_{\gS_i}$ on $\gS_i$ converges to 
$\mu_i|_{\gS_i}$ and $\mu_i|_{\T^n\tim\R^{n+m}\stm \gS_i}=0$. 
Note also that $\mu_i|_{\T^n\tim\R^{n+m}\stm \gS_i}=0$ if and only if 
$\supp \mu_i\subset \gS_i$.  

Thus, we need only to prove that the set $\widetilde\P:=\{\mu|_{\gS}\mid 
\mu\in \P(\gS,a)\}$, where $\mu|_\gS:=(\mu_i|_{\gS_i})_{i\in\I}$, 
is a compact mterizable space with the topology of weak convergence of 
measures.

As noted prior to the lemma, it is enough to prove that $\widetilde\P$ 
is a subset of a closed ball of $\prod_{i\in\I}\M(\gS_i)$ in the norm topology 
and it is closed in the weak convergence of 
measures.   

Recall that the Banach space 
$\prod_{i\in\I}\M(\gS_i)$ has the total variation norm  
$\|\nu\|=\sum_{i\in\I}|\nu_i|(\gS_i)$ for $\nu=(\nu_i)$. 
If $\mu=(\mu_i)\in \M_+(\T^n\tim\R^{n+m})^m$
and $\supp \mu_i\subset \gS_i$ for all $i\in\I$, then  
\[
\du{\mu,\1}=\sum_{i\in\I}\mu_i(\gS_i)=\|\mu|_\gS\|. 
\]
This shows that the closed ball \ $
B:=\{\nu=(\nu_i)\in \prod_{i\in\I}\M(\gS_i)\mid 
\|\nu\|\leq a\}$\ contains $\widetilde\P$.

It remains to show that $\widetilde\P$ is closed in the weak convergence of measures.  For this, as noted at the beginning, we need only to prove that 
$\P(\gS,a)$ is closed in the weak convergence of  measures. 

Now, let $\mu^j=(\mu^j_i)\in \P(\gS,a)$ for all $j\in\N$ and assume that 
the sequence of $\mu^j$ converges weakly in the sense of measures
to $\mu=(\mu_i)\in \M(\T^n\tim\R^{n+m})^m$. 
We already know that $\supp \mu_i\subset \gS_i$ for all $i\in\I$
and that the sequence of $\mu^j|_\gS$ converges to $\mu|_\gS\in B$ 
weakly in the sense of measures. It follows that 
\[
a\geq \|\mu|_{\gS}\|=\sum_{i\in\I}|\mu_i|(\gS_i),
\]  
and, moreover, that for any $i\in\I$ and nonnegative function $\psi\in C_\rmb(\T^n\tim\R^{n+m})$, 
\[
\du{\mu_i,\psi} =\lim_{j\to\infty}\du{\mu^j_i,\psi} \geq 0. 
\]
From these, we see that $\mu\in \M_+(\T\tim\R^{n+m})^m$ and 
$\du{\mu,\1}\leq a$, and conclude that $\mu\in\P(\gS,a)$.
\eproof

\begin{lemma} \label{Plambda} Let $\gl>0$ and $\gS=(\gS_i)_{i\in\I}$ be a 
collection of compact subsets $\gS_i$ of $\T^n\tim\R^n\tim(\R^m\cap Y_i)$. Let $\P^\gl(\gS)$ denote the collection of all $\mu=(\mu_i)\in 
\P^\gl$ such that $\supp \mu_i\subset \gS_i$ 
for every $i\in\I$.      
Then, $\P^\gl(\gS)$ is a compact metrizable space  
with the topology of weak convergence of measures. 
\end{lemma}

\bproof Let $\P(\gS,a)$ denote the set defined in Lemma \ref{PS} for $a>0$. 
For $\mu=(\mu_i)\in\P^\gl(\gS)$, since $\supp \mu_i$ are compact for all $i\in\I$, 
it is clear that $\du{\mu,(|\xi|+|\eta|)\1}<\infty$.  
Since $S^\gl(\eta)\geq \gl$ for all $\eta\in\bigcup_{i\in\I}Y_i$, if  $\mu=(\mu_i)\in\P^\gl(\gS)$, then 
\[
1=\du{\mu,S^\gl\1}
=\sum_{i\in\I}\int_{\gS_i} S^\gl(\eta)\mu_i(dxd\xi d\eta)
\geq \gl \sum_{i\in\I}\mu_i(\gS_i) =\gl \du{\mu, \1},
\]
which implies that $\P^\gl(\gS)\subset \P(\gS,1/\gl)$. 

It remains to prove that $\P^\gl(\gS)$ is a closed subset of 
$\P(\gS,1/\gl)$. 
Let $\mu^j=(\mu^j_i)\in\P^\gl(\gS)$ for $j\in\N$. Assume that 
the sequence $(\mu^j)$ converges weakly in the sense of measures 
to $\mu=(\mu_i)\in\P(\gS,1/\gl)$.  

For the proof of the lemma, we need only to show that $\du{\mu,S^\gl\1}=1$. We easily check that 
\[
\du{\mu,S^\gl\1}=\sum_{i\in\I}\int_{\gS_i}S^\gl(\eta) \mu_i(dxd\xi d\eta)
=\lim_{j\to\infty}
 \sum_{i\in\I}\int_{\gS_i} S^\gl(\eta) \mu^j_i(dx d\xi d\eta)
=\lim_{j\to\infty}\du{\mu^j,S^\gl\1}=1,
\]
which finishes the proof. 
\eproof

It is a consequence of the lemma above that $\P^\gl(K_1,K_2)$ is a compact 
metrizable space with the topology of weak convergence of measures.

\bproof [Proof of Theorem \ref{gpr-thm2}] 
In view of \erf{Two_chara_Pl} and Lemma \ref{gpr-thm1}, it is enough to prove that
\beq \label{gpr4}
v^\gl_k(z)=\min_{\mu\in \tilG\pr(z,k,\gl)}\du{\mu,L}.
\eeq

We intend to show that
\beq\label{gpr5}
\sup_{(\phi,u)\in\tilF(\gl)}\inf_{\nu\in\P^\gl(K_1,K_2)}
\du{\nu,L-\phi+(u_k(z)-v_k^\gl(z))S^\gl\1}= 0.
\eeq

We postpone the proof of \erf{gpr5} and, assuming temporarily that
\erf{gpr5} is valid,  we prove that \erf{gpr4} holds.

To this end, we see easily that $\P^\gl(K_1,K_2)$ is a convex subset of 
a vector space $\M(\T^n\tim\R^{n+m})$ and that, by Lemma \ref{cone}, $\tilF(\gl)$ is a convex subset of $\prod_{i\in\I}C(\T^n\tim Z_i)\tim C(\T^n)^m$.  Observe as well that the functional: 
\[
\P^\gl(K_1,K_2)\ni \nu\mapsto \du{\nu,L-\phi+(u_k(z)-v_k^\gl(z))S^\gl\1}\in\R
\]
is convex and continuous, in the topology of weak convergence of measures for any $(\phi,u)\in\tilF(\gl)$, and the functional: 
\[
\tilF(\gl)\ni (\phi,u)\mapsto \du{\nu,L-\phi+(u_k(z)-v_k^\gl(z))S^\gl\1}\in\R
\]
is concave, as well as continuous, for any $\nu\in\P^\gl(K_1,K_2)$. 

By Lemma \ref{Plambda}, the set $\P^\gl(K_1,K_2)$ is a compact space 
with the topology of weak convergence of measures.    
Hence, we may apply the minimax theorem (Proposition \ref{minimax} or \cite{Te,Si}), 
to deduce from \erf{gpr5} that
\beq \label{gpr7}\bald
0=\sup_{(\phi,u)\in\tilF(\gl)}\min_{\nu\in\P^\gl(K_1,K_2)}&
\du{\nu,L-\phi+(u_k(z)-v_k^\gl(z))S^\gl\1}
\\&=\min_{\nu\in\P^\gl(K_1,K_2)}\sup_{(\phi,u)\in\tilF(\gl)}
\du{\nu,L-\phi+(u_k(z)-v_k^\gl(z))S^\gl\1}.
\eald
\eeq
Observe by using the cone property of $\tilF(\gl)$ that
\[
\sup_{(\phi,u)\in\tilF(\gl)}\du{\nu, u_k(z)S^\gl\1 -\phi}
=\bcases
0 & \IF \ \nu\in \tilG\pr(z,k,\gl), \\[3pt]
\infty&\IF \ \nu\in \P^\gl(K_1,K_2) \stm \tilG\pr(z,k,\gl).
\ecases
\]
This and \erf{gpr7} yield
\[\bald
0&\,=\min_{\nu\in\P^\gl(K_1,K_2)}\sup_{(\phi,u)\in\tilF(\gl)}
\du{\nu,L-\phi+(u_k(z)-v_k^\gl(z))S^\gl\1}
\\&\,
=\min_{\nu\in\tilG\pr(z,k,\gl)}
\du{\nu,L-v_k^\gl(z)S^\gl\1}
=\min_{\nu\in\tilG\pr(z,k,\gl)}
\du{\nu,L}-v_k^\gl(z),
\eald
\]
which proves \erf{gpr4}. 

It remains to show \erf{gpr5}.
Note 
that
\[\bald
\sup_{(\phi,u)\in\tilF(\gl)}\inf_{\nu\in\P^\gl(K_1,K_2)}&
\du{\nu,L-\phi+(u_k(z)-v^\gl_k(z))S^\gl\1}
\\&\geq \inf_{\nu\in \P^\gl(K_1,K_2)}
\du{\nu,L-\phi+(u_k(z)-v_k^\gl(z))S^\gl\1} \Big|_{(\phi,u)=(L,v^\gl)}=0.
\eald\]
Hence, we only need to show that
\beq\label{gpr6}
\sup_{(\phi,u)\in\tilF(\gl)}\inf_{\nu\in \P^\gl(K_1,K_2)}
\du{\nu,L-\phi+(u_k(z)-v_k^\gl(z))S^\gl\1} \leq 0.
\eeq

 For this, we argue by contradiction and thus suppose that \erf{gpr6}
does not hold. Accordingly, we have
\[
\sup_{(\phi,u)\in\tilF(\gl)}\inf_{\nu\in\P^\gl(K_1,K_2)}
\du{\nu,L-\phi+(u_k(z)-v_k^\gl(z))S^\gl\1}> \ep
\]
for some $\ep>0$.  We may select $(\phi,u)\in\tilF(\gl)$ so that
\[
\inf_{\nu\in\P^\gl(K_1,K_2)}\du{\nu,L-\phi+(u_k(z)-v_k^\gl(z))S^\gl\1}
> \ep.
\]
That is, for any $\nu\in\P^\gl(K_1,K_2)$, we have
\[
\du{\nu,L-\phi+(u_k(z)-v_k^\gl(z))S^\gl\1}>\ep=\du{\nu,\ep S^\gl\1}.
\]

According to Lemma \ref{Dirac_measure}, the measure $(S^\gl)^{-1}\gd_{(x,\xi,\eta)} e_i$ is in $\P^\gl(K_1,K_2)$ for every $(x,\xi,\eta)\in\T^n\tim Z_i$ and $i\in\I$.  
Plugging all such $\nu=(S^\gl)^{-1}\gd_{(x,\xi,\eta)} e_i\in \P^\gl(K_1,K_2)$ into the above, we find that
\[
(L_i-\phi_i)(x,\xi,\eta)+(u_k(z)-v_k^\gl(z)-\ep)S^\gl(\eta)>0 \ \ \FORALL (x,\xi,\eta)\in\T^n\tim Z_i,\, i\in\I.
\]
Hence, setting $w:=u-(u_k(z)-v_k^\gl(z)-\ep)\1$,
we have
\[\bald
&\gl w_i(x)+\xi\cdot p+\eta\cdot w(x) -L_i(x,\xi,\eta)
\\&\,=\gl u_i(x) +\xi\cdot p +\eta\cdot u(x)
-(u_k(z)-v_k^\gl(z)-\ep)S^\gl(\eta)
-L_i(x,\xi,\eta)
\\&\,<\gl u_i(x) +\xi\cdot p +\eta\cdot u(x)
-\phi_i(x,\xi,\eta)
\eald
\]
for all $ (x,p,\xi,\eta)\in\T^n\tim\R^n\tim Z_i$ and $ i\in \I$.
This ensures that $w$ is a subsolution of
\[
\gl w+H[w]=0 \ \ \IN \T^n.
\]
By Theorem~\ref{thm2-1}, we get \ $
u(x)-(u_k(z)-v_k^\gl(z)-\ep)\1\leq v^\gl(x)$ for all  $x\in\T^n.$ 
The $k$-th component of the last inequality, evaluated at $x=z$, yields an obvious contradiction, which proves that \erf{gpr6} holds.
\eproof

\def\rmb{\mathrm{b}}

\section{Green-Poisson measures: the general case}

We now remove the hypothesis \erf{compact-continuous} in Theorem \ref{gpr-thm2}
and establish the following theorem.

\begin{theorem}\label{gpg-thm1} Assume \erf{coercive}--\erf{monotone}.
Let $(z,k,\gl)\in \T^n\tim\I\tim (0,\,\infty)$ and $v^\gl\in C(\T^n)^m$ be the solution of \erf{Pl}. 
Then
there exists $\mu\in\frak{C}(z,k,\gl)$ such that
\beq\label{gpg1}
v_k^\gl(z)= \du{\mu,L}=\min_{\nu\in \frak{C}(z,k,\gl)}\du{\nu,L}.
\eeq
\end{theorem}

The 
theorem above guarantees the existence of a Green-Poisson measure associated with any $(z,k,\gl)\in\T^n\tim
\I\tim(0,\,\infty)$.

In what follows we fix a $(z,k,\gl)\in\T^n\tim\I\tim(0,\,\infty)$. 
According to Theorem \ref{thm2-1}, the unique solution of \erf{Pl} is Lipschitz continuous on $\T^n$. With this in mind, we fix a constant $C>0$ and consider the condition that
\beq\label{gpg3}
|v^\gl(x)|+|Dv^\gl(x)|\leq C \ \ \text{ a.e. }x\in \T^n.
\eeq
We choose a function $h\in C^1(\R^{n}\tim\R^m)$ so that
\beq\label{gpg4}
\left\{\bald
&\text{$h$ is nonnegative and convex on } \ \R^{n}\tim\R^m,
\\
& h(p,u)=0 \ \ \ \text{ if and only if }\ \  |p|+|u|\leq C,
\\
&
\lim_{|p|+|u|\to\infty}(|p|+|u|)^{-1}h(p,u)=\infty.
\eald\right.
\eeq
Also, we choose a compact convex set $Q\subset \R^{n+m}$
such that for all $(x,i,p,u)\in\T^n\tim\I\tim\R^{n+m}$,
\beq\label{gpg5}
\pl_{(p,u)}H_i(x,p,u)\subset Q \ \ \IF |p|+|u|\leq C,
\eeq
where $\pl_{(p,u)}H_i$ denotes the subdifferential
of the convex function: $(p,u)\mapsto H_i(x,p,u)$.

\begin{theorem} \label{gpg-thm2} Assume \erf{coercive}--\erf{monotone}. Let $v^\gl$ be the solution of \erf{Pl}
and assume that \erf{gpg3} is satisfied for some
constant $C>0$. Let $Q$ be a compact convex subset of
$\R^{n+m}$ such that \erf{gpg5} holds.
Assume that there exists $\mu\in\frak{C}(z,k,\gl)$
such that
\[ 
v^\gl_k(z)=\du{\mu,L}.
\] 
Then
\[
\supp \mu_i \subset \T^n\tim [Q\cap(\R^n\tim Y_i)] \ \ \FOR i\in\I.
\]
\end{theorem}

We recall some basic properties related to the subdifferentials
of $H$ and $L$.

\begin{lemma}\label{basic_subdiff} Assume \erf{convex}.
Let $(x,i)\in\T^n\tim\I$. 

\noindent 
\emph{(i)} We have
\[
\pl_{(p,u)}H_i(x,p,u)\not=\emptyset \ \ \FOR (p,u)\in\R^n\tim\R^m.
\]
\emph{(ii)} Let $(p,u),\,(\xi,\eta)\in\R^{n+m}$. The following three statements are equivalent each other.
\begin{enumerate}
\item[\emph{(a)}] $(\xi,\eta)\in\pl_{(p,u)}H_i(x,p,u)$.
\item[\emph{(b)}] $(p,u)\in\pl_{(\xi,\eta)}L_i(x,\xi,\eta)$.
\item[\emph{(c)}] $H_i(x,p,u)+L_i(x,\xi,\eta)=\xi\cdot p+\eta\cdot u$.
\end{enumerate}
\end{lemma}

\bproof (i) \ Since $(p,u)\mapsto H_i(x,p,u)$ is continuous and convex in $\R^{n+m}$,
it is locally Lipschitz continuous (see \cite[Theorem B.3]{Ishort}) and hence almost everywhere
differentiable (see \cite[Theorem F.1]{Ishort})
in $\R^{n+m}$. Fix any $(p,u)\in\R^{n+m}$ and choose a sequence
of points $(p^k,u^k)\in\R^{n+m}$ converging to $(p,u)$ such that
$(p,u)\mapsto H_i(x,p,u)$ is differentiable at $(p^k,u^k)$ for all $k\in\N$.
Set $(\xi^k,\eta^k)=D_{p,u}H_i(x,p^k,u^k)$ for $k\in\N$.
The local Lipschitz continuity of $H_i(x,\cdot,\cdot)$ allows us to assume that
$(\xi^k,\eta^k)_{k\in\N}$ is bounded and, moreover, convergent to some
$(\xi^0,\eta^0)\in\R^{n+m}$ after passing to a subsequence. Since
\[
H_i(x,p^k+q,u^k+r)\geq H_i(x,p^k,u^k)+\xi^k\cdot q+\eta^k\cdot r \ \ \FOR (q,r)\in\R^{n+m},
\]
sending $k\to\infty$ yields
\[
H_i(x,p+q,u+r)\geq H_i(x,p,u)+\xi^0\cdot q+\eta^0\cdot r \ \ \FOR (q,r)\in\R^{n+m},
\]
which shows that $(\xi^0,\eta^0)\in\pl_{(p,u)}H_i(x,p,u)$ and $\pl_{(p,u)}H_i(x,p,u)
\not=\emptyset$.

We here skip to prove (ii) and leave it to the reader to consult
\cite[Theorem 23.5]{Ro}
or \cite[Theorem B.2]{Ishort}.
\eproof

\begin{lemma} \label{C to G} Assume \erf{coercive}--\erf{monotone}. Let 
$\gl\in[0,\,\infty)$, and let $u\in\Lip(\T^n)$ be a subsolution of \erf{Pl}. 
If $\gl>0$, then $\,u_k(z)\leq \du{\mu,L}\,$ for all $(z,k)\in\T^n\tim\I$ and $\mu\in\fC(z,k,\gl)$, 
and, if $\gl=0$, then $\,0\leq \du{\mu,L}\,$ for all $\mu\in\fC(0)$. 
\end{lemma}

We remark that, in the above, $\du{\mu,L}$ can be $+\infty$. 

The proof below is almost identical to the first part of the proof of Lemma \ref{gpr-thm1}. 

\bproof
In view of the continuity and convex property of $H$, mollifying $u$, we may choose, for each $\ep>0$,  
a function $u^\ep\in C^1(\T^n)^m$ such that 
$\gl u^\ep+H[u^\ep]\leq \ep\,\1$ in $\T^n$ and $\|u-u^\ep\|_\infty<\ep$.  Hence, we have
\[
\gl u^\ep+\xi\cdot Du^\ep(x)+\eta\cdot u^\ep(x)\,\1 \leq L(x,\xi,\eta)+\ep\,\1.
\]
When $\gl>0$, fixing $(z,k)\in\T^n\tim\I$, recalling \erf{gpr1}, and integrating the inequality above with respect to $\mu\in\fC(z,k,\gl)$, we obtain
\[\bald
u^\ep_k(z)&\,=\du{\mu,\, \xi\cdot Du^\ep+\eta\cdot u^\ep\,\1+\gl u^\ep}
\leq \du{\mu,\, L+\ep\,\1}.
\eald
\]
Similarly, if $\gl=0$, then we get for any $\mu\in\fC(0)$, 
\[
0=\du{\mu, \xi\cdot Du^\ep+\eta\cdot u^\ep} 
\leq \du{\mu, L+\ep\,\1}. 
\] 
Taking the limit as $\ep\to 0$,  we finish the proof. 
\eproof 

 In the proof below, an essential step is to construct a new Hamiltonian, say, $\widetilde H$ satisfying \erf{coercive}--\erf{monotone} 
such that $v^\gl$ is a solution of \erf{Pl}, with 
$H$ replaced by $\widetilde H$, and such that, if $|p|+|u|> C$, then $\widetilde H(x,p,u)>H(x,p,u)$. Notice that the function $H(x,p,u)+h(p,u)$, with $h$ satisfying \erf{gpg4}, on $\T^n\tim\R^{n+m}$ 
does not satisfy the monotonicity \erf{monotone}. 

\bproof[Proof of Theorem \ref{gpg-thm2}]
 Let $h\in C^1(\R^n\tim\R^m)$ be a function having the properties in \erf{gpg4}. 
We set $G^h(x,p,u)=H(x,p,u)+h(p,u)$ for
$(x,p,u)\in\T^n\tim\R^n\tim\R^m$. 
Let $K^h=(K^h_i)_{i\in\I}$ be
the Lagrangian of $G^h$, and, since $G^h$ grows superlinearly as
$|p|+|u| \to \infty$, we see that $K^h\in C(\T^n\tim\R^n\tim\R^m)^m$.  Note that
\[
G^h\geq H \ \ \AND \ \ K^h\leq L \ \ \ \ON
\T^n\tim \R^n\tim \R^m.
\]

According to Proposition \ref{prop2.1}, $G^h$ does not satisfy \erf{monotone}, and we  
need to modify $G^h$, to remove the drawback.  
We note by Proposition \ref{prop2.1} that 
\[
L(x,\xi,\eta)+\0_{Y}(\eta)=L(x,\xi,\eta)
\ \ \FOR (x,\xi,\eta)\in\T^n\tim\R^n\tim\R^m,
\]
where $\0_Y:=(\0_{Y_i})_{i\in\I}$. Hence, we have 
\[
L^h(x,\xi,\eta):=K^h(x,\xi,\eta)+\0_{Y}(\eta)\leq
L(x,\xi,\eta)
\ \ \FOR (x,\xi,\eta)\in\T^n\tim\R^n\tim\R^m.
\]
Let $H^h=(H^h_i)_{i\in\I}$ be the Hamiltonian
of $L^h$, and note that 
\[ 
H\leq H^h\leq G^h \ \ \ON \ \T^n\tim\R^n\tim\R^m.
\] 
In particular, we have
\[
H(x,p,u)=H^h(x,p,u)=G^h(x,p,u) \ \ \IF |p|+|u|\leq C,
\]
which shows, together with \erf{gpg3}, that $v^\gl$ is a solution
of $\gl u+H^h[u]=0$ in $\T^n$. It is clear that
$H^h$ satisfies \erf{coercive} and \erf{convex}.
Moreover, $H^h$ satisfies \erf{monotone} due to
Proposition \ref{prop2.1}.

Now, since $L\geq L^h$ on $\T^n\tim\R^n\tim\R^m$, it follows immediately that
\[
v_k^\gl(z)=\du{\mu,L}\geq\du{\mu,L^h}.
\]
Since $\gl v^\gl+H^h[v^\gl]=0$ in $\T^n$, thanks to Lemma \ref{C to G}, we get
\[
v_k^\gl(z)\leq \du{\mu,L^h}.
\]
Combining these yields
\[
v_k^\gl(z)=\du{\mu,L}=\du{\mu,L^h}.
\]
Consequently, we have 
\beq\label{=&ineq}
\du{\mu,L-L^h}=0 \ \ \AND \ \ L\geq L^h.
\eeq

Noting that for all $i\in\I$, $Y_i$ is a closed subset of $\R^m$ and $L_i=\infty$ on $\T^n\tim\R^n\tim(\R^m\stm Y_i)$, and $L_i-L^h_i$ is lower semicontinuous on
$\T^n\tim\R^n\tim Y_i$, we easily deduce from \erf{=&ineq} that 
\[
\supp \mu_i\subset \{(x,\xi,\eta)\in \T^n\tim\R^n\tim Y_i\mid L_i(x,\xi,\eta)=L_i^h(x,\xi,\eta)\}.
\]

It remains to show that for all $i\in\I$,
\beq\label{gpg8}
 \{(x,\xi,\eta)\in \T^n\tim\R^n\tim Y_i\mid L_i(x,\xi,\eta)=L^h_i(x,\xi,\eta)\}\subset \T^n\tim Q.
\eeq
To do this, we fix $i\in\I$ and
\[(x,\xi,\eta)\in \T^n\tim\R^n\tim Y_i
\]
such that $L_i(x,\xi,\eta)=L_i^h(x,\xi,\eta)$, set
$\gz=(\xi,\eta)$ and show that $\gz\in Q$.
We argue by contradiction and thus suppose that
$\gz\not\in Q$.

Note that, since $\gz\in \R^n\tim Y_i$,
\beq\label{gpg9}
K_i^h(x,\gz)=L_i^h(x,\gz)=L_i(x,\gz).
\eeq
In view of Lemma \ref{basic_subdiff}, (i) applied to
$K^h$, we can
select $q_\gz=(p_\gz,u_\gz)\in \pl_{(\xi,\eta)}K^h_i(x,\gz)$, which implies
by the convex duality (Lemma \ref{basic_subdiff}, (ii))  that $\gz\in\pl_{(p,u)}G^h_i(x,q_\gz)$ and
\beq \label{gpg10}
K^h_i(x,\gz)+G_i^h(x,q_\gz)=\gz\cdot q_\gz.
\eeq

We claim that $h(q_\gz)>0$. Indeed, if, to the contrary, $h(q_\gz)=0$, then we have
$|p_\gz|+|u_\gz|\leq C$ by \erf{gpg4} and,
by \erf{gpg5}, \erf{gpg9}, and \erf{gpg10},
\[
\pl_{(p,u)}H_i(x,q_\gz)\subset Q \ \ \AND \ \ 
\gz\cdot q_\gz=K^h_i(x,\gz)+G_i^h(x,q_\gz)
=L_i(x,\gz)+H_i(x,q_\gz),
\]
which imply by Lemma \ref{basic_subdiff}, (ii) that
\[\gz\in\pl_{(p,u)}H_i(x,q_\gz) \subset Q.
\]
This contradicts the choice of $\gz$, which confirms that $h(q_\gz)>0$.

Now, we observe that
\[\bald
L_i(x,\gz)&\,\geq \gz\cdot q_\gz -H_i(x,q_\gz)
=\gz\cdot q_\gz -G^h_i(x,q_\gz)+h(q_\gz)
\\&\,
=K_i^h(x,\gz)+h(q_\gz)>K_i^h(x,\gz)=L_i(x,\gz),
\eald\]
which is a contradiction, and we conclude that \erf{gpg8} is valid. The proof is complete.
\eproof

\def\tilH{\widetilde H} 
\def\tilL{\widetilde L}
\def\fC{\frak{C}}

In the following proof of Theorem \ref{gpg-thm1}, we approximate the Hamiltonian $H(x,p,u)$ by Hamiltonians which satisfy \erf{coercive}--\erf{monotone} and \erf{compact-continuous}. 
In the first step of the approximation of $H$, we follow 
the argument in the proof above, with $h$ replaced by $h/r$, with $r\in\N$. 
In the proof above, the function $G^h(x,p,u)$ has the superlinear growth
in $(p,u)$ because of the addition of $h$ and its nice effect is the continuity 
of $K^h$ on $\T^n\tim\R^{n+m}$. The continuity on $\T^n\tim\R^{n+m}$ of the Lagrangians of the approximating Hamiltonians, obtained in the first step, is important for the second and final step of building the approximating Hamiltonians, which have at most the linear growth due to \erf{compact-continuous}.    

\bproof[Proof of Theorem \ref{gpg-thm1}] We choose
a constant $C>0$ and a compact convex set $Q\subset \R^{n+m} $ so that \erf{gpg3} and \erf{gpg5} hold. We may assume that $Q = Q_1\tim Q_2$ for 
some $Q_1\subset  \R^n$ and $Q_2\subset\R^m$, where,
moreover, $Q_1$ is a neighborhood of the origin of $\R^n$.
Let $h\in C^1(\R^{n+m})$ be a function
satisfying \erf{gpg4}.
As in the proof of Theorem~\ref{gpg-thm2}, we define
sequences $(H^r)_{r\in\N},\,(L^r)_{r\in\N},\,
(G^r)_{r\in\N},\,(K^r)_{r\in\N}$ of functions, with
$h$ replaced by $h/r$. That is, $G^r=(G^r_i)_{i\in\I}$ is defined by
\[
G^r_i(x,p,u)=H_i(x,p,u)+\fr 1 r \, h(p,u) \ \ \FOR (x,p,u)
\in\T^n\tim\R^{n+m},
\]
$K^r$ is the Lagrangian of 
$G^r$, $L^r$  is given by
\[
L^r(x,\xi,\eta)=K^r(x,\xi,\eta)+\0_{Y}(\eta) \ \ \FOR (x,\xi,\eta)
\in\T^n\tim\R^{n+m},
\]
and $H^r$ is the Hamiltonian of 
$L^r$.
We have already checked  in the proof of Theorem~\ref{gpg-thm2} that
$H^r$ satisfies \erf{coercive}--\erf{monotone},
$v^\gl$ is a solution of $\gl v^\gl+H^r[v^\gl]=0$
in $\T^n$, and $L^r\in\prod_{i\in\I}C(\T^n\tim\R^n\tim Y_i)$.
Moreover, it is easily seen that for
$(x,p,u)\in\T^n\tim\R^{n+m}$ and $i\in\I$, if $|p|+|u|\leq C$,
\beq\label{gpg11}
H(x,p,u)=H^r(x,p,u)=G^r(x,p,u) \ \ \AND \ \ \pl_{(p,u)}H^r_i(x,p,u)\subset Q.
\eeq

Next we define function $H^r_Q=(H^r_{Q,i})_{i\in\I}$
as the Hamiltonian of the function
\[
L^r_Q(x,\xi,\eta):=L^r(x,\xi,\eta)+\0_{Q}(\xi,\eta).
\]
Note by Lemma \ref{basic_subdiff}, (ii) that for
$(x,i,p,u)\in\T^n\tim\I\tim\R^{n+m}$ and $\gz\in\R^{n+m}$,  if
\[
\gz\in\pl_{(p,u)}H_{Q,i}^r(x,p,u),
\]
then
\[
(p,u)\in\pl_{(\xi,\eta)}L_{Q,i}^r(x,\gz),
\]
and hence, by the definition of $L_{Q,i}^r$, we have \ $
\gz\in Q.$ 
That is, we have
\[
\pl_{(p,u)}H^r_{Q,i}(x,p,u)\subset Q\ \ \FOR
(x,i,p,u)\in\T^n\tim\I\tim\R^{n+m}.
\]
It is now easy to see that $H^r_Q$ satisfies
\erf{coercive}, \erf{convex} and \erf{compact-continuous}.
Note also by the inclusion in \erf{gpg11} that
if $|p|+|u|\leq C$,
\[\bald
H^r_i(x,p,u)&\,=\max_{(\xi,\eta)\in Q}(p\cdot\xi+u\cdot \eta -L^r_i(x,\xi,\eta))
\\&\,=\max_{(\xi,\eta)\in (\R^n\tim Y_i)\cap  Q}(p\cdot\xi+u\cdot \eta -L^r_i(x,\xi,\eta))
=H^r_{Q,i}(x,p,u).
\eald
\]

We may now invoke Theorem \ref{gpr-thm2}, to conclude that there is $\mu\in\fC(z,k,\gl)$
such that
\beq\label{gpg12}
v_k^\gl(z)=\du{\mu,L^r_Q}=\min_{\nu\in\fC(z,k,\gl)}
\du{\nu,L^r_Q}.
\eeq
Theorem~\ref{gpg-thm2} and \erf{gpg11} ensure that
for any minimizer $\nu=(\nu_i)_{i\in\I}\in\fC(z,k,\gl)$ of  the
optimization in \erf{gpg12}, we have the property 
\[
\supp \nu_i \subset \T^n\tim[(\R^n\tim Y_i)\cap Q].
\]

For each $r\in\N$, we select a minimizer $\mu^r=(\mu_i^r)\in\fC(z,k,\gl)$ of  the
optimization in \erf{gpg12}. Since $\supp\mu_i^r\subset \T^n\tim Q$, we see 
immediately from\erf{gpg12} that 
\beq \label{gpg12+}
v_k^\gl(z)=\du{\mu^r,L^r}.
\eeq

In view of Lemma \ref{Plambda}, we may assume that $(\mu^r)_{r\in\N}$, 
after passing to a subsequence which is denoted again by the same symbol, converges weakly in the sense of measures to a measure 
$\mu=(\mu_i)\in\P^\gl$ having the property that $\supp \mu_i \subset \T^n\tim[(\R^n\tim Y_i)\cap Q]$ for all $i\in\I$.  

The weak convergence of $(\mu^r)$ implies that 
\[\left\{\bald
&\du{\mu,S^\gl\1}=1,  
\\&\psi_k(z)=\du{\mu, \xi\cdot D\psi_i+\eta\cdot \psi\1+\gl \psi}
 \ \ \FORALL \psi=(\psi_i)\in C^1(\T^m)^m.
\eald\right.
\]
These ensure that $\mu\in \fC(z,k,\gl)$.

It is easily checked that, as $r\to\infty$, $\,
K^r(x,\xi,\eta) \to L(x,\xi,\eta)\,$
monotonically pointwise. Since $K^r \leq K^{r+1}$ and $K^r\leq L^r$ 
for $r\in\N$,
we obtain from \erf{gpg12+},  
\[ 
v_k^\gl(z)\geq \du{\mu^r,K^q}\geq \du{\mu^r, j\wedge K^q} \ \ \IF r\geq q,\ \ 
\FORALL j,q\in\N,
\] 
where $j\wedge K^q:=\left(\min\{j,\,K_i^q\}\right)_{i\in\I}\in C_\rmb(\T^m\tim\R^{n+m})^m$. 
Sending $r\to\infty$ yields 
\[ 
v_k^\gl(z)\geq \du{\mu, j\wedge K^q} \ \ \FORALL j, q\in\N.
\] 
By the monotone convergence theorem, after sending $j,q\to\infty$, we obtain
\[
v_k^\gl(z)\geq \du{\mu,L},
\]
while, by Lemma \ref{C to G}, we have 
\[
v_k^{\gl}(z)\leq \inf_{\nu\in\fC(z,k,\gl)}\du{\nu,L}. 
\]
Thus, we conclude that 
\[
v_k^\gl(z)=\du{\mu,L}=\min_{\nu\in\fC(z,k,\gl)}\du{\nu,L}. \qedhere
\]
\eproof

\section{A convergence result for the vanishing discount problem}

We study the asymptotic behavior of the solution $v^\gl$ of \erf{Pl}, with $\gl>0$,
as $\gl\to 0$.

\begin{theorem} \label{conv-thm1}
Assume \erf{coercive}--\erf{ergodic}.
Let $v^\gl$ be the solution of \erf{Pl} for $\gl>0$. Then there exists a solution $v^0$ of \emph{(P$_0$)}
such that  the functions $v^\gl$ converge to $v^0$ in $C(\T^n)^m$ as
$\gl\to 0+$.
\end{theorem}

\begin{lemma}\label{conv-thm2} Under the hypotheses of Theorem \ref{conv-thm1},
there exists a constant $C_0>0$ such that for any $\gl>0$,
\beq\label{conv1}
|v_i^\gl(x)|\leq C_0 \ \ \FOR (x,i)\in\T^n\tim\I.
\eeq
\end{lemma}

\bproof Let $v_0=(v_{0,i})_{i\in\I}\in \Lip(\T^n)^m$ be a solution of (P$_0$). Choose
a constant $C_1>0$ so that
\[
|v_{0,i}(x)|\leq C_1 \ \ \FOR (x,i)\in\T^n\tim\I,
\]
and observe by Lemma \ref{lem1} that the functions $v_0+C_1\1$
and $v_0-C_1\1$ are a supersolution and a subsolution of (P$_0$), respectively.
Noting that $v_0+C_1\1\geq 0$ and $v_0-C_1\1\leq 0$, we deduce that
$v_0+C_1\1\geq 0$ and $v_0-C_1\1\leq 0$ are a supersolution and a subsolution of
\erf{Pl}, respectively, for any $\gl>0$. By comparison (Theorem \ref{thm2-1}),
we see that, for any $\gl>0$,
$v_0-C_1\1\leq v^\gl\leq v_0+C_1\1$ on $\T^n$ and, moreover,
$-2C_1\1\leq v^\gl\leq 2C_1\1$ on $\T^n$. Thus, \erf{conv1} holds with
$C_0=2C_1$.
\eproof

\begin{lemma}\label{conv-thm3} Under the hypotheses of Theorem \ref{conv-thm1},
the family $(v^\gl)_{\gl\in(0,\,1)}$ is equi-Lipschitz continuous on $\T^n$.
\end{lemma}

\bproof According to Lemma \ref{conv-thm2}, we may choose a constant $C_0>0$ so that
\[
|v^\gl_i(x)|\leq C_0 \ \ \FOR (x,i,\gl)\in\T^n\tim\I\tim (0,\,\infty).
\]
Hence, as $v^\gl$ is a solution of \erf{Pl},   we deduce by \erf{coercive} that there
exists a constant $C_1>0$ such that the
functions $v_i^\gl$, with $\gl\in(0,\,1)$, are subsolutions of
$|Du|\leq C_1$ in $\T^n$. As is well-known, this implies that the $v_i^\gl$ are
Lipschitz continuous on $\T^n$ with $C_1$ as their Lipschitz bound.
\eproof

We remark that one can show, with a slightly more elaboration,
the equi-Lipschitz property of $(v^\gl)_{\gl>0}$ in the above lemma.

\begin{theorem}\label{conv-thm4} Let $(z,k)\in\T^n\tim\I$. Assume \erf{coercive}--\erf{ergodic}.
For any $\gl>0$, let $v^{\gl}$ be the solution of \erf{Pl} and 
$\mu^\gl \in\fC(z,k,\gl)$ a minimizer in \erf{gpg1}. Then, for any sequence $(\gl_j)_{j\in\N}$ of positive numbers converging to zero, there exists a subsequence of $(\gl_j)$, which is denoted again by the same symbol, such that, as $j\to\infty$,
\[
\gl_j \mu^{\gl_j} \to \nu^0
\]
weakly in the sense of measures for some $\nu^0=(\nu^0_i)_{i\in\I}\in\fC(0)$,
and $\nu^0$ satisfies
\beq\label{conv2}
0=\du{\nu^0,L}=\min_{\nu\in\fC(0)}\du{\nu,L}.
\eeq
\end{theorem}

We call any minimizing measure $\nu^0\in\fC(0)$ in \erf{conv2} a \emph{Mather measure.} The set of all Mather measures $\nu^0\in\fC(0)$ is denoted by
$\frak{M}(L)$.  See, for example, \cite{Mat, Man, DFIZ, DZ2} for some work related to 
Mather measures.  
Notice that the limit measure $\nu^0$ in Theorem \ref{conv-thm4}
is a Mather measure.
It should be noted that, in our formulation, the existence of a Mather measure is trivial since $0\in\fC(0)$.

\bproof We fix $(z,k)\in \T^n\tim\I$. By Theorem \ref{gpg-thm1}, for each $\gl>0$
there exists $\mu^\gl=(\mu_i^\gl)_{i\in\I}\in \fC(z,k,\gl)$ such that
\beq\label{conv3}
\gl v_k^\gl(z)= \du{\gl \mu^\gl,L}.
\eeq

By Lemmas \ref{conv-thm2} and \ref{conv-thm3},
there is a constant $C>0$ such that for any $\gl\in(0,\,1)$,
\[
|v^\gl(z)|+|Dv^\gl(x)|\leq C \ \ \text{ a.e. }
x\in \T^n.
\]
We choose
a closed ball $Q
\subset\R^{n+m}$ so that
\erf{gpg5} holds with $C$ given above. Thanks to Theorem~\ref{gpg-thm2},
we find that
\[
\supp \mu_i^\gl\subset \T^n\tim[Q\cap (\R^n\tim Y_i)]\ \ \FOR (i,\gl)\in\I\tim(0,\,1).
\]

Noting that $S^\gl(\eta)\geq\gl$ for $\eta\in Y_i$ and $\supp \mu_i^\gl \subset \T^n\tim\R^n\tim Y_i$ for all $i\in\I$, we observe that
\[
\du{\gl\mu^\gl,\1}=\du{\mu^\gl, \gl \1} \leq \du{\mu^\gl,S^\gl \1}=1.
\]
Hence, by applying Lemma \ref{PS} and passing to a subsequence, we may assume that the sequence $(\gl_j \mu^{\gl_j})_{j\in\N}\subset \M_+(\T^n\tim\R^{n+m})^m $  converges weakly in the sense of measures to some  $\nu^0=(\nu^0_i)\in 
\M_+(\T^n\tim\R^{n+m})^m$ having properties 
$\du{\nu^0,\1}\leq 1$ and $\supp \nu_i^0 \subset \T^n\tim[Q\cap(\R^n\tim Y_i)]$ 
for all $i\in\I$. It is an immediate consequence that $\nu^0\in\P^0$. 

Since $\mu^\gl\in\fC(z,k,\gl)$, we have 
\[
\gl \psi_k(z)=\du{\gl \mu^\gl, \xi\cdot D\psi+\eta\cdot \psi\1+\gl \psi} 
\ \ \FORALL \psi=(\psi_i)\in C^1(\T^n)^m.
\]
Sending $\gl \to 0$ along the sequence $(\gl_j)$ yields 
\[
0=\du{\nu^0, \xi\cdot D\psi+\eta\cdot \psi\1} 
\ \ \FORALL \psi=(\psi_i)\in C^1(\T^n)^m,
\]
which concludes that $\nu^0\in\fC(0)$. 

For any  function
$\phi\in C_\rmb(\T^n\tim\R^{n+m})^m$ such that $\phi\leq L$, we see from \erf{conv3} that
\[
0\geq \du{\nu^0,\phi}.
\]
Moreover, by approximating $L$ monotonically  from below by bounded continuous functions and applying  the monotone convergence theorem, we deduce that
\[ 
0\geq \du{\nu^0,L}.
\] 
By Lemma \ref{C to G}, we have $\,0\leq\du{\nu,L}$ for all $\nu\in \fC(0)$. 
It is now clear that \erf{conv2} holds. 
\eproof

Let $\cV$ denote the set of accumulation points $v
=(v_i)_{i\in\I}\in C(\T^n)^m$ of $(v^\gl)_{\gl>0}$ in the space $C(\T^n)^m$ as $\gl\to 0$.
Note by the stability of the viscosity property under uniform
convergence that any $v\in\cV$ is a solution of (P$_0$).
Let $\cW$ denote the set of those solutions $w\in C(\T^n)^m$ of (P$_0$)
which satisfy
\beq \label{conv5}
\du{\nu,w}\leq 0 \ \ \FORALL \ \nu\in\fM(L).
\eeq

\bproof[Proof of Theorem \ref{conv-thm1}]
In view of the Ascoli-Arzela theorem, Lemmas \ref{conv-thm2} and \ref{conv-thm3} assure that the family $(v^\gl)_{\gl\in(0,\,1)}$ is relatively compact in $C(\T^n)^m$.  In particular, the set $\cV$ is nonempty.

If $\cV$ is a singleton, then it is obvious that the whole family
$(v^\gl)_{\gl>0}$ converges to the unique element of $\cV$ in $C(\T^n)^m$
as $\gl\to 0$.

We need only to show that $\cV$ is a singleton.  For this, we first show that
\beq\label{conv6}
\cV\subset\cW.
\eeq
To see this, let $v\in\cV$ and $\nu\in\fM(L)$.
Choose a sequence $(\gl_j)_{j\in\N}$ of positive
numbers converging to zero such that $(v^{\gl_j})_{j\in\N}$ converges to $v$ in
$C(\T^n)^m$.
Since $(L-\gl v^\gl,v^\gl)\in\cF(0)$, $\nu\in\fC(0)$, 
and $\du{\nu,L}=0$, 
using Lemma \ref{C to G},  
we get 
\[
0\leq \du{\nu,L-\gl v^\gl}=
\du{\nu,L}-\du{\nu, \gl v^\gl}
=-\gl\du{\nu,v^\gl},
\]
which yields, after dividing  by $\gl>0$ and then sending $\gl \to 0$
along $\gl=\gl_j$,
\[
\du{\nu,v}\leq 0.
\]
This proves \erf{conv5}, which ensures the inclusion \erf{conv6}.

 Next, we show that
\beq\label{conv7}
w\leq v \ \ \FORALL\ w\in\cW,\,v\in\cV.
\eeq

To check this, it is enough to show that
for any $v\in\cV$, $w\in\cW$ and $(z,k)\in\T^n\tim\I$, the inequality $w_k(z)\leq v_k(z)$ holds.

Fix any $v\in\cV$ and $w\in\cW$ and $(z,k)\in\T^n\tim\I$.
Select a sequence $(\gl_j)_{j\in\N}\subset(0,\,\infty)$ converging to zero
so that
\[
v^{\gl_j} \to v \ \ \IN C(\T^n)^m \ \ \hb{ as } j\to\infty.
\]
By Theorem \ref{gpg-thm1}, there exists a sequence $(\mu^j)_{j\in\N}$  such that   for $j\in\N$,
\beq\label{conv8}
\mu^j\in \fC(z,k,\gl_j) \ \ \AND \ \
v_k^{\gl_j}(z)=\du{\mu^j, L}.
\eeq
In view of Theorem \ref{conv-thm4}, we may assume by passing to a subsequence if necessary that, as $j\to\infty$,
\[
\gl_j \mu^j \to \nu \ \ \hb{ weakly in the sense of measures}
\]
for some $\nu=(\nu_i)_{i\in\I}\in\fM(L)$.

Now, note that $(L+\gl_j w, w)\in\cF(\gl_j)$ and infer by Lemma \ref{C to G} and \erf{conv8} that
\[
w_k(z)\leq \du{\mu^j, L+\gl_j w}
=v_k^{\gl_j}(z)+\gl_j\du{\mu^j, w}.
\]
Sending $j\to\infty$ now yields
\[
w_k(z)\leq v_k(z)+\du{\nu,w}.
\]
This together with \erf{conv5} shows that $w_k(z)\leq v_k(z)$, which
ensures that \erf{conv7} holds.
Noting that \erf{conv7} combined with \erf{conv6}  shows that $w\leq v$ for all $v,w\in\cV$, that is, $\cV$ is a singleton. The proof is complete.  \eproof

Reviewing the proof above, we conclude easily the following proposition, which is a
generalization of \cite[Theorem 3.8]{DFIZ} (see also \cite[Proof of Theorem 1]{DZ2}).

\begin{corollary} \label{conv-thm5} Under the assumptions and notation of Theorem \ref{conv-thm1}, the limit function $v^0=(v^0_i)_{i\in\I}$ can be represented as
\[
v_i^0(x)=\max\{w_i(x)\mid w=(w_i)\in\cW\} \ \ \FOR\ x\in\T^n.
\]
\end{corollary}

The proof of Corollary \ref{conv-thm5}, with $\cW$ replaced by
\[
\cW^-=\{w\in C(\T^n)^m\mid w \text{ is a subsolution and satisfies \erf{conv5}}\},
\]
shows also that, under the hypotheses and notation of Corollary \ref{conv-thm5},
\[
v^0_i(x)=\max\{w_i(x)\mid w=(w_i)\in\cW^-\} \ \ \FOR x\in\T^n.
\]

\section{Ergodic problem}

Remark that, given a Hamiltonian $H$, condition \erf{ergodic} is not satisfied in general.  We consider the problem
of finding an $m$-vector $c=(c_i)_{i\in\I}\in\R^m$ and a function
$u=(u_i)_{i\in\I}\in C(\T^n)^m$ such that $u$ is a solution of the $m$-system
\beq\label{erg1}
H[u]=c \ \ \IN \T^n,
\eeq
which is stated componentwise as
\[
H_i(x,Du_i(x),u(x))=c_i \ \ \IN \T^n\ \FOR i\in\I.
\]
We call this problem the \emph{ergodic problem} for $H$.

If the ergodic problem has a solution $c\in\R^m$ and $u\in C(\T^n)^m$, then we may apply the main convergence result (Theorem \ref{conv-thm1}) to \erf{Pl}, with $H$ replaced by $H_c:=H-c$.  As noted in the introduction, this change of Hamiltonians, in general, does not help analyze the vanishing discount problem for the original system \erf{Pl}.  

However, if $H$ satisfies a certain additional condition, then the argument of switching from the Hamiltonian $H$ to $H_c$ makes sense for the vanishing 
discount problem for \erf{Pl}. For instance, given a 
solution $(c,u)\in \R^m\tim C(\T^n)^m$ of \erf{erg1}, assume that  
the equality 
\beq\label{homo-c}
H(x,p,v+tc)=H(x,p,v) 
\eeq
holds for all $t\in\R$ and $(x,p,v)\in\T^n\tim\R^{n+m}$. It is easily seen that
if $v^\gl\in C(\T^n)^m$ is a solution of \erf{Pl}, then $w^\gl:=v^\gl+\gl^{-1}c$ 
is a solution of \erf{Pl}, with $H_c$ in place of $H$.  This is a situation where one 
can apply Theorem \ref{conv-thm1}, to observe the convergence of $v^\gl+\gl^{-1}c$
as $\gl\to 0$.

In the next result, we do not need the convexity or monotonicity of $H$, and we assume only \erf{coercive}.

For $R>0$ and $r>0$, we set
\[\bald
\ga_R(r)&\,=\inf\{H_i(x,p,u)\mid (x,i)\in\T^n\tim\I,\ u\in B_R^m,\ p\in \R^n\stm B_r^n\},
\\ \gb_R&\,=\sup\{H_i(x,0,u)\mid (x,i)\in\T^n\tim \I,\ u\in B^m_R   \}.
\eald
\]
The constants $\ga_R(r)$ and $\gb_R$ are finite by the continuity of $H_i$ and \erf{coercive}.
It is clear that for any $R>0$, the function $r\mapsto \ga_R(r)$ is nondecreasing in $(0,\,\infty)$ and  diverges to infinity as $r\to \infty$.

\begin{theorem} \label{thm-erg} Assume \erf{coercive} and that there
exists a constant $R>0$ such that
\beq \label{erg2}
\gb_R<\ga_R\left(\textstyle \fr {2R}{\sqrt n}\right).
\eeq
Then problem \erf{erg1} has a solution $(c,u)\in\R^m\tim C(\T^n)^m$.
\end{theorem}

We remark that a result 
similar to the above has been established in \cite[Theorem 1.2]{JMT} in the case of a scalar Hamilton-Jacobi equation.

Roughly speaking, the condition \erf{erg2} in the theorem above is satisfied 
for a large $R>0$ if the growth of $H_i(x,p,u)$ in $p$ is higher in a certain sense  than that in $u$ 
as $|(p,u)|\to \infty$.   
In the case of linear coupling (and hence, $H_i$ have the form of \erf{linear}), 
it is obvious that if for all $i\in\I$, the functions $G_i(x,p)$ have the superlinear  growth, i.e., 
satisfy 
\[
\lim_{R\to\infty}\inf_{(x,p)\in \T^n\tim(\R^n\stm B_R^n)}\fr{G_i(x,p)}{|p|}=\infty,
\]
then condition \erf{erg2} is valid. 
We refer to \cite[Theorem 2.12]{DZ1}, 
\cite[Theorem 17]{Icoupling} for results, in the linear coupling case, 
similar to but more subtle than the theorem above. 

Condition \erf{erg2} is also valid when, as a direct generalization of the linear coupling case, $H_i, i\in\I$ have superlinear growth in $p$, i.e., for any $r>0$,
\[
\lim_{R\to\infty}\inf_{(x,p,u)\in \T^n\tim(\R^n\stm B_R^n)\tim B_r^m}\fr{H_{i}(x,p,u)}{|p|}=\infty,
\]
and uniformly Lipschitz dependence in $u$, i.e., there exists $\Theta>0$ such that for any $u,v\in\R^m$
\[
|H_{i}(x,p,u)-H_{i}(x,p,v)|\leq\Theta|u-v|.
\]

\bproof  We choose $R>0$ so that \erf{erg2} holds and select $\gl>0$ so that
\beq \label{erg3}
\gb_R+\gl R<\ga_R\left(\textstyle \fr{2R}{\sqrt n}\right).
\eeq

Let $u\in C(\T^n)^m$ and consider the uncoupled $m$-system for $v=(v_i)_{i\in\I}$:
\beq \label{erg4}
\gl(v_i(x)-u_i(x))+H_i(x,Dv_i(x),u(x))=0 \ \ \IN \T^n \ \FOR i\in\I.
\eeq
The functions $(x,p)\mapsto H_i(x,p,u(x))$ are continuous and coercive and, hence,
the standard theory of viscosity solutions (also, Theorem \ref{thm2-1} applied to each single equations) guarantees that
\erf{erg4} has a unique solution $v=(v_i)_{i\in\I}$ and the functions $v_i$
are Lipschitz continuous on $\T^n$.

For any $u\in C(\T^n)^m$, let $v=(v_i)_{i\in\I}\in C(\T^n)^m$ be the solution of  \erf{erg4}. We set
\[
Tu :=v-\min_{\T^n}v,
\]
where
\[
\min_{T^n}v:=(\min_{x\in\T^n} v_i(x))_{i\in\I} \in\R^m,
\]
which gives a mapping $T$ from $C(\T^n)^m$ to $C(\T^n)^m$. Because of the stability
of viscosity solutions under the uniform convergence
and the uniqueness of solution of \erf{erg4}, we easily deduce that $T$ is a
continuous mapping on the Banach space $C(\T^n)^m$, with norm $\|u\|_\infty :=\max_{x\in\T^n}|u(x)|$.

Now, fix $u$ so that
\[
\|u\|_\infty\leq R \ \ \AND \ \ u(x)\geq 0 \ \ \FORALL x\in\T^n,
\]
and observe that the function $w(x):=-\gl^{-1}\gb_R \1$ is a subsolution of \erf{erg4}. Indeed,
we have
\[
\gl(w_i(x)-u_i(x)) +H_i(x,Dw_i(x),u(x)) \leq -\gb_R+H_i(x,0,u(x))\leq 0 \ \ \FOR (x,i)\in\T^n\tim\I.
\]
By the standard comparison theorem, we have
\[
-\fr{\gb_R}{\gl}\1\leq v.
\]

Noting that $v_i$ is Lipschitz continuous and hence it is almost everywhere differentiable,
we compute at any point $x$ of differentiability of $v_i$ that,  if $Dv_i(x)\not=0$,
\[
0\geq \gl\left(-\fr{\gb_R}{\gl}-u_i(x)\right)+\ga_R(|Dv_i(x)|)
\geq -\gb_R-\gl R+\ga_R(|Dv_i(x)|),
\]
and observe by the choice of $\gl$ that,  if $Dv_i(x)\not=0$,
\[
\ga_R(|Dv_i(x)|)\leq \gb_R+\gl R<\ga_R\left(\textstyle\fr {2R}{\sqrt  n}\right),
\]
which yields
\[
|Dv_i(x)|<\fr {2R}{\sqrt n} \ \ \text{ a.e. in }\T^n\ \FOR i\in\I,
\]
and moreover
\[
0\leq v_i(x)-\min_{T^n} v_i \leq R \ \ \text{ for all }(x,i)\in\T^n\tim\I.
\]
Thus, we conclude that
\[\bald
&\|D(Tu)_i\|_{L^\infty(\T^n)}:=\mathop{\operatorname{ess\, sup}}_{\T^n}|D(Tu)_i|\leq R \ \FOR i\in\I,
\\& (Tu)(x)\geq  0 \ \ \FOR x\in\T^n \quad \AND \quad \|Tu\|_\infty\leq R.
\eald
\]

We set
\[
K=\{u\in C(\T^n)^m\mid u\geq 0,\ \|u\|_\infty\leq R,\ \|Du_i\|_{L^\infty(\T^n)}\leq R \ \FORALL i\in\I\},
\]
and note that $K$ is a compact convex subset of $C(\T^n)^m$.  The above observations
show that $T$ maps $K$ into $K$. The Schauder fixed point theorem guarantees that there
is a fixed point $u\in K$ of $T$. Let $v$ be the a solution of \erf{erg4}, with the fixed point $u$.
By the definition of $T$, we have
\[
u=Tu=v-\min_{\T^n}v,
\]
and $u$ solves
\[
\gl \min_{\T^n} v_i +H_i(x,Du_i,u)=0 \ \ \IN \T^n \ \ \FOR i\in\I.
\]
That is, the pair $(-\gl \min_{\T^n}v, \,u)$ is a solution of \erf{erg2}.
\eproof

\section*{Acknowledgement}
The authors would like to thank Wenjia Jing of Tsinghua University, 
who kindly shared the ideas for the proof of \cite[Theorem 1.2]{JMT} before its
publication,  which was a great help for them to establish Theorem 16.

The authors would like to thank the anonymous referees for their careful reading 
of and critical and useful comments on the original version of this paper, which 
have helped significantly to improve the presentation.

HI was supported in part by the JSPS Grants KAKENHI  No. 16H03948 and 
No. 18H00833 and by the NSF Grant No. 1440140 while in residence at the Mathematical Sciences Research Institute in Berkeley, California, in October 2018.
HI thanks the Department of Mathematics at the Sapienza University of Rome for financial support and its hospitality while his visit there in May 6--June 5  2019. 
LJ was supported by the National Natural Science Foundation of China (Grant No. 11901293 and No. 11571166) and Start-up Foundation of Nanjing University of Science and Technology (No. AE89991/114).

\bye

\bib{DZ1}{article}{
   author={Davini, Andrea},
   author={Zavidovique, Maxime},
   title={Aubry sets for weakly coupled systems of Hamilton-Jacobi
   equations},
   journal={SIAM J. Math. Anal.},
   volume={46},
   date={2014},
   number={5},
   pages={3361--3389},
   issn={0036-1410},
   review={\MR{3265180}},
   doi={10.1137/120899960},
}

\bib{MT1}{article}{
   author={Mitake, Hiroyoshi},
   author={Tran, Hung V.},
   title={Remarks on the large time behavior of viscosity solutions of
   quasi-monotone weakly coupled systems of Hamilton-Jacobi equations},
   journal={Asymptot. Anal.},
   volume={77},
   date={2012},
   number={1-2},
   pages={43--70},
   issn={0921-7134},
   review={\MR{2952714}},
}

\bib{MT2}{article}{
   author={Mitake, H.},
   author={Tran, H. V.},
   title={A dynamical approach to the large-time behavior of solutions to
   weakly coupled systems of Hamilton-Jacobi equations},
   language={English, with English and French summaries},
   journal={J. Math. Pures Appl. (9)},
   volume={101},
   date={2014},
   number={1},
   pages={76--93},
   issn={0021-7824},
   review={\MR{3133425}},
   doi={10.1016/j.matpur.2013.05.004},
}

\bib{MT3}{article}{
   author={Mitake, Hiroyoshi},
   author={Tran, Hung V.},
   title={Homogenization of weakly coupled systems of Hamilton-Jacobi
   equations with fast switching rates},
   journal={Arch. Ration. Mech. Anal.},
   volume={211},
   date={2014},
   number={3},
   pages={733--769},
   issn={0003-9527},
   review={\MR{3158806}},
   doi={10.1007/s00205-013-0685-x},
}

\bib{Var}{book}{
   author={Varga, Richard S.},
   title={Matrix iterative analysis},
   series={Springer Series in Computational Mathematics},
   volume={27},
   edition={Second revised and expanded edition},
   publisher={Springer-Verlag, Berlin},
   date={2000},
   pages={x+358},
   isbn={3-540-66321-5},
   review={\MR{1753713}},
   doi={10.1007/978-3-642-05156-2},
}

\bib{CLLN}{article}{
   author={Camilli, Fabio},
   author={Ley, Olivier},
   author={Loreti, Paola},
   author={Nguyen, Vinh Duc},
   title={Large time behavior of weakly coupled systems of first-order
   Hamilton-Jacobi equations},
   journal={NoDEA Nonlinear Differential Equations Appl.},
   volume={19},
   date={2012},
   number={6},
   pages={719--749},
   issn={1021-9722},
   review={\MR{2996426}},
   doi={10.1007/s00030-011-0149-7},
}

\bib{CGT}{article}{
   author={Cagnetti, Filippo},
   author={Gomes, Diogo},
   author={Tran, Hung Vinh},
   title={Adjoint methods for obstacle problems and weakly coupled systems
   of PDE},
   journal={ESAIM Control Optim. Calc. Var.},
   volume={19},
   date={2013},
   number={3},
   pages={754--779},
   issn={1292-8119},
   review={\MR{3092361}},
   doi={10.1051/cocv/2012032},
}

\bib{FrH1}{article}{
   author={Freidlin, Mark},
   author={Hu, Wenqing},
   title={Smoluchowski-Kramers approximation in the case of variable
   friction},
   note={Problems in mathematical analysis. No. 61},
   journal={J. Math. Sci. (N.Y.)},
   volume={179},
   date={2011},
   number={1},
   pages={184--207},
   issn={1072-3374},
   doi={10.1007/s10958-011-0589-y},
}

\bib{FrHW}{article}{
   author={Freidlin, Mark},
   author={Hu, Wenqing},
   author={Wentzell, Alexander},
   title={Small mass asymptotic for the motion with vanishing friction},
   journal={Stochastic Process. Appl.},
   volume={123},
   date={2013},
   number={1},
   pages={45--75},
   issn={0304-4149},
   doi={10.1016/j.spa.2012.08.013},
}

\bib{GiTr}{book}{
   author={Gilbarg, David},
   author={Trudinger, Neil S.},
   title={Elliptic partial differential equations of second order},
   series={Classics in Mathematics},
   note={Reprint of the 1998 edition},
   publisher={Springer-Verlag, Berlin},
   date={2001},
   pages={xiv+517},
   isbn={3-540-41160-7},
}

\bib{E1}{article}{
   author={Evans, Lawrence C.},
   title={The perturbed test function method for viscosity solutions of
   nonlinear PDE},
   journal={Proc. Roy. Soc. Edinburgh Sect. A},
   volume={111},
   date={1989},
   number={3-4},
   pages={359--375},
   issn={0308-2105},
   doi={10.1017/S0308210500018631},
}
\bib{E2}{article}{
   author={Evans, Lawrence C.},
   title={Periodic homogenisation of certain fully nonlinear partial
   differential equations},
   journal={Proc. Roy. Soc. Edinburgh Sect. A},
   volume={120},
   date={1992},
   number={3-4},
   pages={245--265},
   issn={0308-2105},
   doi={10.1017/S0308210500032121},
}

\bib{Li83}{article}{
   author={Lieberman, Gary M.},
   title={The conormal derivative problem for elliptic equations of
   variational type},
   journal={J. Differential Equations},
   volume={49},
   date={1983},
   number={2},
   pages={218--257},
   issn={0022-0396},
   doi={10.1016/0022-0396(83)90013-X},
}
\bib{Patrizi}{article}{
   author={Patrizi, Stefania},
   title={Principal eigenvalues for Isaacs operators with Neumann boundary
   conditions},
   journal={NoDEA Nonlinear Differential Equations Appl.},
   volume={16},
   date={2009},
   number={1},
   pages={79--107},
   issn={1021-9722},
   doi={10.1007/s00030-008-7042-z},
}
